\documentclass{amsart}

\usepackage{amssymb}

\def\nin{\noindent}
\def\bc{\begin{center}}
\def\ec{\end{center}}
\def\be{\begin{equation}}
\def\ben{\begin{enumerate}}
\def\een{\end{enumerate}}
\def\bfg{\begin{figure}}
\def\efg{\end{figure}}
\def\bq{\begin{quote}}
\def\eq{\end{quote}}
\def\bd{\begin{description}}
\def\ed{\end{description}}
\def\this{i.\ e.\ } %that is

\def\h{\hbar}
\def\p{\partial}
\def\w{\wedge}

\def\dim{\operatorname{dim}}

\newcommand{\CC}{{\mathbb C}}

\newcommand{\ZZ}{{\mathbb Z}}
\newcommand{\QQ}{{\mathbb Q}}

\newcommand{\lan}{\langle}
\newcommand{\ran}{\rangle}
\newcommand{\gL}{\Lambda}

\newcommand{\gS}{\Sigma}

\newcommand{\calo}{{\mathcal O}}

\newcommand{\calt}{{\mathcal T}}

\newcommand{\M}{\overline{\mathcal M}}

\newcommand{\ev}{\operatorname{ev}}
\newcommand{\ft}{\operatorname{ft}}
\renewcommand{\QQ}{\mathbb Q}
\renewcommand{\CC}{\mathbb C}
\renewcommand{\t}{\mathbf t}
\newcommand{\q}{\mathbf q}
\newcommand{\f}{\mathbf f}
\newcommand{\g}{\mathbf g}

\newcommand{\D}{\mathcal D}

\newcommand{\A}{\mathcal A}
\renewcommand{\a}{\alpha}
\renewcommand{\b}{\beta}
\renewcommand{\c}{\gamma}
\renewcommand{\d}{\delta}

\renewcommand{\H}{{\mathcal H}}
\newcommand{\F}{{\mathcal F}}
\newcommand{\G}{{\mathcal G}}

\newcommand{\C}{{\mathcal C}}
\renewcommand{\L}{{\mathcal L}}
\renewcommand{\O}{{\Omega}}
\newcommand{\1}{{\bf 1}}
\newcommand{\s}{{\mathbf s}}
\renewcommand{\S}{{\Sigma}}
\newcommand{\EE}{{\mathbf E}}
\renewcommand{\l}{{\lambda}}
\newcommand{\ch}{\operatorname{ch}}
\newcommand{\cc}{{\mathbf c}}
\newcommand{\e}{{\mathbf e}}

\newcommand{\0}{{\mathbf 0}}

\newcommand{\str}{\operatorname{str}}
\newcommand{\sdet}{\operatorname{sdet}}
\newcommand{\hh }{\hat{\ }} 

\renewcommand{\(}{\left(}
\renewcommand{\)}{\right)}
\newcommand{\Td}{\operatorname{Td}}
\newcommand{\Tdv}{\Td^{\vee}}
\usepackage{diagrams}

\newcommand{\<}{\left<}
\renewcommand{\>}{\right>}

\begin{document}

\title[Quantum Riemann -- Roch, Lefschetz and Serre ]
{Quantum Riemann -- Roch, Lefschetz and Serre}
 
\author{Tom Coates \and Alexander  Givental} 
\address{UC Berkeley} 
\thanks{Research is partially supported by NSF Grant DMS-0072658} 

%\date{October 12, 2001}

\begin{abstract}

Given a holomorphic vector bundle $E:EX\to X$ over a compact K\"ahler manifold,
one introduces twisted GW-invariants of $X$ replacing virtual fundamental 
cycles of moduli spaces of stable maps $f: \Sigma \to X$ by their cap-product 
with a chosen multiplicative characteristic class of 
$H^0(\Sigma, f^* E) - H^1(\Sigma, f^*E)$. 
Using the formalism \cite{Gi0} of quantized quadratic hamiltonians, 
we express the descendent potential for the twisted theory in terms of that 
for $X$. The result (Theorem $1$) is a consequence of Mumford's Riemann -- 
Roch -- Grothendieck formula \cite{Mu,FP} applied to the universal stable map. 

When $E$ is concave, and the inverse $\CC^{\times}$-equivariant Euler class is 
chosen, the twisted theory yields GW-invariants of $EX$. The ``non-linear 
Serre duality principle'' \cite{Gi1, Gi2} expresses GW-invariants of $EX$ via 
those of the supermanifold $\Pi E^*X$, where the Euler class and $E^*$ replace
the inverse Euler class and $E$. We derive from Theorem $1$ the nonlinear Serre
duality in a very general form (Corollary $2$). 

When the bundle $E$ is convex, and a submanifold $Y\subset X$ is defined by 
a global section, the genus $0$ GW-invariants of $\Pi E X$ coincide with 
those of $Y$. We prove a ``quantum Lefschetz hyperplane section principle'' 
(Theorem $2$) expressing genus $0$ GW-invariants of a complete intersection $Y$
via those of $X$. This extends earlier results \cite{BCKS, Kim, Ber, Lee, Gath}
and yields most of the known mirror formulas for toric complete intersections. 
\end{abstract}

\maketitle

\nin {\bf Introduction.} The mirror formula of Candelas {\em et al}
\cite{COGP} for the virtual numbers $n_d$ of degree $d=1,2,3,...$
holomorphic spheres on a quintic $3$-fold $Y\subset X=\CC P^4$ can be
stated \cite{Gi} as the coincidence of the $2$-dimensional
cones over the following two curves in $H^{even}(Y;\QQ )= \QQ [P]/(P^4)$:
\[ J_{Y}(\tau):= e^{P\tau}+\frac{P^2}{5}\sum_{d>0}n_dd^3 
\sum_{k>0}\frac{e^{(P+kd)\tau}}{(P+kd)^2} \]
and
\[ I_Y(t) = \sum_{d\geq 0} e^{(P+d)t}\frac{(5P+1)(5P+2)...(5P+5d)}
{(P+1)^5(P+2)^5 ... (P+d)^5} .\]    
The new proof given in this paper shares with earlier work 
\cite{Gi1, Kim, LLY, Ber, Lee, Gath} the formulation of sphere 
counting in a hypersurface $Y \subset X$ as a problem in Gromov -- Witten
theory of $X$. \\

Gromov -- Witten invariants of a compact almost K\"ahler manifold $X$
are defined as intersection numbers in moduli spaces $X_{g,n,d}$ 
of stable pseudo-holomorphic maps $f:\Sigma \to X$. 
All results of this paper can be stated and hold true 
in this generality, but we prefer to stay on the firmer ground of 
algebraic geometry where the most applications belong.  

Given a holomorphic vector bundle $E$ over a complex projective manifolds $X$ and 
an invertible multiplicative characteristic class $\cc $ of complex vector bundles,
we introduce {\em twisted} Gromov -- Witten invariants as intersection
indices in $X_{g,n,d}$ with the characteristic classes
$\cc (E_{g,n,d})$ of the virtual bundles 
$E_{g,n,d}= ``H^0(\Sigma, f^*E)\ominus H^1(\Sigma, f^*E) ''$.
The ``quantum Riemann -- Roch theorem'' (Theorem $1$) expresses the 
twisted Gromov -- Witten invariants (of any genus) and their gravitational
descendents via untwisted ones.

The totality of gravitational descendents 
in the genus $0$ Gromov -- Witten theory of $X$ can be encoded by a semi-infinite 
cone $\L_{X}$ in the cohomology algebra of $X$ with coefficients in the field
of Laurent series in $1/z$ (see Section $6$). 
Another such cone corresponds to each twisted theory. Let $\L_{E}$ be
the cone corresponding to the total Chern class
\[ \cc = \l^{\dim }+c_1\l^{\dim -1}+...+c_{\dim} .\]
Theorem $1$ specialized to this case 
says that the cones $\L_X$ and $\L_{E}$ are related 
by a linear transformation. It is described in terms of the  
stationary phase asymptotics $a_{\rho} (z)$ of the oscillating integral
\[  \frac{1}{\sqrt{2\pi z}} \int_0^{\infty} e^{\frac{-x+(\l+\rho) \ln x}
{z}} dx \]  
as multiplication in the cohomology algebra by $\prod_i a_{\rho_i}(z)$, where 
$\rho_i $ are the Chern roots of $E$.

Assuming $E$ to be a line bundle, we derive a 
``quantum hyperplane section theorem'' (Theorem $2$). It is more
general than the earlier versions \cite{BCKS, Kim, Lee, Gath} in the
sense that the restrictions $t\in H^{\leq 2}(X;\QQ)$ on the space of parameters
and $c_1(E)\leq c_1(X)$ on the Fano index are removed. 

In the quintic case when $X=\CC P^4$ and $\rho =5P$, the cone 
$\L_X$ is known to contain the curve
\[ J_X(t)=\sum_{d\geq 0} \frac{e^{(P+zd)t/z}}{(P+z)^5...(P+zd)^5}, \]
and Theorem $2$ says that the cone $\L_{E}$ contains the curve 
\[ I_{E}(t) = \sum_{d\geq 0} e^{(P+zd)t/z}
\frac{(\l+5P+z) ... (\l+5P+5dz)}{(P+z)^5...(P+dz)^5}. \]
One obtains the quintic mirror formula
by passing to the limit $\l=0$. \\

The idea of deriving mirror formulas by 
applying the Grothendieck -- Riemann -- Roch theorem to universal
stable maps is not new. Apparently this was the initial plan
of M. Kontsevich back in $1993$. In $2000$, we had a chance to
discuss a similar proposal with R. Pandharipande. We would like
to thank these authors as well as  A. Barnard and 
A. Knutson for helpful conversations.

The second author is grateful to D. van Straten for the 
invitation to the workshop ``Algebraic aspects of mirror symmetry''
held at Kaiserslautern in June $2001$. The discussions at the workshop 
and particularly the lectures on ``Variations of semi-infinite Hodge 
structures'' by S. Barannikov proved to be very useful in our work 
on this project. 

\bigskip

\nin {\bf 1. Generating functions.} 
Let $X$ be a compact complex projective manifold of complex dimension $D$.
Denote by $X_{g,n,d}$ the moduli orbispace of genus $g$, $n$-pointed stable maps 
\cite{Kn1, BM} to $X$ of degree $d$, where $d \in H_2(X;\ZZ)$.
The moduli space is compact and can be equipped \cite{BF, LT, S}  
with a (rational-coefficient) {\em virtual fundamental cycle} 
$[X_{g,n,d}]$ of complex dimension $n+(1-g)(D-3)+\int_{d} c_1(TX)$.

The {\em total descendent potential} of $X$ is a generating function
for Gromov-Witten invariants. It is defined as
\begin{equation} {\D}_X := \exp \( \sum \h^{g-1} {\F}^g_X \) ,
\label{D} \end{equation}
where ${\F}^g_X$ is the {\em genus $g$ descendent potential},
\begin{equation} \F^g_X=\sum_{n,d} \frac{Q^d}{n!} \int_{[X_{g,n,d}]} 
 (\sum_{k=0}^{\infty} (\ev_1^*t_k) \psi_1^k ) ...  
(\sum_{k=0}^{\infty} (\ev_n^*t_k) \psi_n^k ) .
\label{F^g} \end{equation}
Here $\psi_i$ is the $1$-st Chern class of the universal 
cotangent line bundle over $X_{g,n,d}$ corresponding to the $i$-th marked 
point, $\ev_i^*t_k$ are the pull-backs by the evaluation map 
$\ev_i: X_{g,n,d}\to X$ at the $i$-th marked point of the cohomology classes
$t_0,t_1,... \in H^*(X,\QQ)$, and $Q^d$ is the representative of $d$ in 
the semigroup ring of degrees of holomorphic curves in $X$.\\

Let $E: EX \to X$ be a holomorphic vector bundle. We regard it as an
element of the Grothendieck group $K^0(X)$. Given $E$, 
one associates to it an element $E_{g,n,d}$ in the Grothendieck group $K^0(X_{g,n,d})$
of coherent orbisheaves as follows. 
Consider the {\em universal stable map}
with the indices $(g,n,d)$ 
\begin{diagram}
X_{g,n+1,d} & \rTo^{\ev} & X \\
\dTo^{\ft} \\
X_{g,n,d} 
\end{diagram}
formed by the maps of forgetting and evaluation at the last marked point.
We pull-back $E$ to the universal curve and then apply the
$K$-theoretic push-forward to $X_{g,n,d}$. Namely, there exists a 
complex $0 \to E^0_{g,n,d} \to E_{g,n,d}^1 \to 0$ of locally free orbisheaves on 
$X_{g,n,d}$ whose cohomology sheaves are respectively $ R^0\ft_* (\ev^* E)$ and  
$R^1 \ft_* (\ev^* E)$. Moreover, the difference 
\[ E_{g,n,d} := [E^0_{g,n,d}] - [E^1_{g,n,d}]   \]
in the Grothendieck group of vector orbibundles does not depend on the choice of the
complex. These facts are based on some standard general results \cite{SGA} about
{\em local complete intersection morphisms}. We refer to Appendix $1$ for a further
discussion of these properties of the maps $\ft$.

A rational invertible multiplicative characteristic 
class of complex vector bundles takes on the form
\begin{equation} \label{cc} \cc (\cdot ) = \exp \{ \sum_{k=0}^{\infty} s_k \ch_k(\cdot) \} 
\end{equation}   
where $\ch_k$ are components of the Chern character, and
$\s = (s_0, s_1, s_2, ...)$ are arbitrary coefficients or indeterminates.
Given such a class and a vector bundle $E\in K^0(X)$, one introduces the
{\em $(\cc ,E)$-twisted descendent potentials} 
$\D_{\cc ,E}$ and $\F^g_{\cc ,E}$ by replacing the virtual 
fundamental cycles $[X_{g,n,d}]$ in (\ref{D}, \ref{F^g}) with the cap-products 
$\cc (E_{g,n,d})\cap [X_{g,n,d}]$. 
For example, the Poincar\'e intersection pairing arises in Gromov -- 
Witten theory as the intersection index in $X_{0,3,0}=X$, and in
the twisted theory therefore takes on the form  
\be \label{ip} (a,b)_{\cc(E)} := \int_{[X_{0,3,0}]}\cc (E_{0,3,0})
\ev_1^*(a) \ev_2^*(1) \ev_3^*(b) = \int_{X} \cc (E)\ a\ b .
\end{equation} 

The genus $0$ potentials $\F_{\cc, E}^0$ when reduced modulo
$Q$ have zero $2$-jet at $\t =0$ since 
$X_{0,n,0}=X\times \M_{0,n}$ and  
$\dim \M_{0,n}=n-3$.
This shows that the twisted potentials $\D_{\cc , E}$ are well-defined
(despite of the occurrence of both $\h$ and $\h^{-1}$ in the exponent of 
(\ref{D})) as formal power series in $\t/\h, Q/\h$ and $\h$. \\
 
We will often assume that all vector bundles
carry the $S^1$-action defined as the fiberwise multiplication by the unitary
scalars. In this case $\ch_k$   
are understood as $S^1$-equivariant characteristic classes 
and all GW-invariants take values in the coefficient ring of 
$S^1$-equivariant cohomology theory. 
We will always identify the ring $H^*(BS^1;\QQ)$ with $\QQ [ \l ]$ where 
$\l $ is the $1$-st Chern class of the Hopf bundle over $\CC P^{\infty}$.    

\bigskip

\nin {\bf 2. Quantization formalism.} 
Theorem $1$ below expresses $\D_{\cc ,E}$ via $\D_X$ in terms of the
formalism of quadratic hamiltonians and their quantization in the Fock space
(see \cite{Gi0}).

Consider the cohomology space $H=H^*(X;\QQ )$ as a super-space equipped
with the non-degenerate symmetric bilinear form defined by the Poincar\'e 
intersection pairing
$( a, b ) = \int_X a b $. Let $\H = H[z,z^{-1}]$ denote the super-space 
of Laurent polynomials in one even indeterminate $z$ with coefficients in $H$.
We equip $\H$ with the even symplectic form
\begin{equation} \label{O} 
 \O (\f,\g) =\frac{1}{2\pi i} \oint ( \f(-z), \g(z)) dz = - 
(-1)^{\bar{\f}\bar{\g}} \O (\g, \f). \end{equation}    
The polarization $\H=\H_{+}\oplus \H_{-}$ defined by the 
Lagrangian subspaces $\H_{+}=H [z]$, $\H_{-}=z^{-1} H [z^{-1}]$ identifies
$(\H, \O)$ with the cotangent bundle $T^*\H_{+}$. Then the standard
quantization convention associates to quadratic hamiltonians $G$ on $(\H, \O)$ 
differential operators $\hat{G}$ of order $\leq 2$ acting on functions on 
$\H_{+}$. 

More precisely, let $\{ q_{\a} \}$ be a $\ZZ_2$-graded coordinate 
system on $\H_{+}$ and $\{ p_{\a} \}$ be the dual coordinate system on 
$\H_{-}$ so that the symplectic structure in these coordinates assumes
the Darboux form
\[  \O (\f,\g) = \sum_{\a} 
[ p_{\a}(\f) q_{\a}(\g) - (-1)^{\bar{p}_{\a}\bar{q}_{\a}} q_{\a}(\f) 
p_{\a}(\g) ]. \] 
For example, when $H$ is the standard one-dimensional Euclidean space then
$\f= \sum q_k z^k + \sum p_k (-z)^{-1-k}$ is such a coordinate system.

In a Darboux coordinate system the quantization convention reads 
\[ (q_{\a}q_{\b})\hh := \frac{q_{\a}q_{\b}}{\h} ,\ \ 
(q_{\a}p_{\b})\hh := q_{\a} \frac{\p}{\p q_{\b}}, \ \ 
(p_{\a}p_{\b})\hh := \h \frac{\p^2}{\p q_{\a}\p q_{\b}} .\]
The quantization is only a {\em projective} representation of the Lie algebra 
of quadratic hamiltonians on $\H$ to the Lie algebra of differential operators.
For quadratic hamiltonians $F$ and $G$ we have
  \[ \{ F, G\}\hh = [ \hat{F}, \hat{G}] + \C (F,G) \]
where $\{ \cdot , \cdot \}$ is the Poisson bracket, $[\cdot , \cdot ]$ is the
super-commutator, and $\C $ is the cocycle
\[
\begin{array}{lcll} 
\C (p_{\a}p_{\b}, q_{\a}q_{\b}) &=& (-1)^{\bar{q}_{\a}\bar{p}_{\b}} \ & 
\text{if}\ \a \neq \b,\\  \C (p_{\a}^2, q_{\a}^2) &=& 
1+(-1)^{\bar{q}_{\a}\bar{p}_{\a}}, & \\ 
\C =0 & & \text{on any other pair} & \text{of quadratic Darboux monomials.}  
\end{array} \]  
We associate the quadratic hamiltonian $(T\f,\f)/2$ to an infinitesimal symplectic 
transformation $T$. If $A, B$ are self-adjoint operators on $H$, then 
the operators $\f \mapsto (A/z) \f$ and $\f \mapsto (Bz) \f $ in $\H$ 
are infinitesimal symplectic transformations, and 
\[ \C (A/z , Bz) =  \str (AB)/2 . \]

\medskip
 
The differential operators act on formal functions 
(with coefficients depending on $\h^{\pm 1}$) 
on the space $\H_{+}$ of vector-polynomials 
$\q = q_0+q_1z+q_2z^2+...$  with the coefficients $q_0,q_1,q_2 ... \in H$.
We will often refer to such functions as {\em elements of the Fock space}. 
We will assume that the ground field 
$\QQ $ of constants is extended
to the Novikov ring $\QQ [[ Q ]]$ (or to $\QQ (\lambda) [[Q]]$ in the
$S^1$-equivariant setting) and  will denote the ground ring by
$\Lambda $. 
      
On the other hand, the potentials $\F^g_X$ are naturally defined as 
formal functions on the 
space of vector-polynomials $\t (\psi) = t_0+t_1\psi+t_2\psi^2+...$ where 
$t_0,t_1,t_2,... \in H$ are cohomology classes of $X$ with coefficients in 
$\gL$, and $\psi$ is an indeterminate to be substituted by  
successive $\psi_i$'s in the definition (\ref{F^g}). We identify $\psi$ with 
$z$, and the total descendent potential (\ref{D}) --- 
with an element of the Fock space by means of the {\em dilaton shift} 
\[ \q (z) = \t (z) - z\ . \] 

The twisted descendent potentials $\D_{\cc , E}$ can be similarly considered 
as elements of the Fock spaces corresponding to the super-space $H$ equipped 
with the twisted inner products (\ref{ip}).
Alternatively, we can identify the inner product spaces 
$(H, (\cdot ,\cdot )_{\cc (E)} )$ with $(H, (\cdot, \cdot ))$
by means of the maps $a \mapsto a \sqrt{\cc (E)}$. Using the corresponding
identification of the Fock spaces we consider
the twisted descendent potentials $\D_{\cc , E}$
as elements of the original Fock space via the convention:
\be \label{ds} \q (z) = \sqrt{\cc (E)} (\t (z) - z) \ .\end{equation}
We obtain therefore a family $\D_{\s} := \D_{\cc , E}$ of 
elements of the Fock space depending on the parameters 
$\s = (s_0,s_1,s_2,...)$. It is easy to see that the potentials $\D_{\cc , E}$ 
are well-defined at least as {\em formal functions of variables $\t$
and parameters $\s$}.
We should stress however that, due to the dilaton shift, 
$\D_{\s}$ as an element of the Fock space is a formal function of $\q$ near 
the shifted origin $\q(z)=-\sqrt{\cc(E)} z$ (drifting therefore with $\s $). 

\bigskip
       
\nin {\bf 3. Quantum Riemann--Roch.} 
The operator of multiplication by $\ch_l(E)$ in the cohomology 
{\em algebra} $H$ of $X$ is self-adjoint with respect to the Poincare pairing. 
Consequently the operator $A$ of multiplication by $\ch_l(E) z^{2m-1}$ in the
algebra $\H = H [z,z^{-1}]$ is anti-self-adjoint with respect to the
symplectic structure $\O$. Thus this operator defines an infinitesimal
symplectic transformation, and $(\ch_l(E) z^{2m-1})\hh $ denotes
the operator on the Fock space defined by quantization of the corresponding 
quadratic hamiltonian $\O (Af, f)/2$. 
  
Let $B_{2m}$ denote  Bernoulli numbers:
$  \frac{x}{1-e^{-x}} =  \frac{x}{2} + \sum_{m\geq 0} 
B_{2m} \frac{x^{2m}}{(2m)!}\ $.

\medskip
  
{\bf Theorem 1.} 
\be \label{T1} 
\exp \{-\frac{1}{24} \sum_{l>0} s_{l-1}\int_X \ch_l(E) c_{D-1}(T_X) \}
\ (\operatorname{sdet} \sqrt{\cc (E)})^{-\frac{1}{24}}\ \D_{\s}  = \end{equation}
\[ \exp \{ \sum_{m>0}\sum_{l\geq 0} s_{2m-1+l} 
 \frac{B_{2m}}{(2m)!} \widehat{( \ch_l(E) z^{2m-1} )} \} \ \ 
\exp \{ \sum_{l>0} s_{l-1} \widehat{(\ch_l (E)/z)}\}\ \D_{\0} .\]

\medskip

{\em Remarks.} (1) The variable $s_0$ is present on the RHS of (\ref{T1}) 
only in the form
$\exp (s_0 \rho /z)\hh $ where $\rho = \ch_1(E)$. For any 
$\rho \in H^2(X)$ the operator $(\rho /z)\hh $ is in fact a
{\em divisor} operator, that is the total descendent potential
satisfies the following {\em divisor equation}:
\be \label{de} (\frac{\rho}{z})\hh \ \D_{\0} = 
\sum \rho_i Q_i\frac{\p}{\p Q_i} \ \D_{\0} 
-\frac{1}{24} \int_X \rho\ c_{D-1}(T_X) \ \D_{\0}.\end{equation}
Here $Q_i$ are generators in the Novikov ring corresponding to a choice of
a basis in $H_2(X)$, and $\rho_i$ are coordinates of $\rho $ in the dual basis.
For $\rho = \ch_1(E)$ the $c_{D-1}$-term cancels with the $s_0$-term on the
LHS of (\ref{T1}). Thus the action of the $s_0$-flow reduces  
to the change $Q^d \mapsto Q^d \exp (s_0 \int_d \rho )$ 
in the descendent potential $\D_{\0}$ combined with the multiplication by the factor
$\exp \{ s_0 (\dim E) /48 \}$ which comes from the super-determinant.     

(2) When $E=\CC $, we have $E_{g,n,d}=\CC - \EE^*_g$ where $\EE_g$ is the
{\em Hodge bundle}. The Hodge bundles are known to
satisfy $\ch_k(\EE_g)=-\ch_k(\EE_g^*)$ (in fact we reprove this in the 
next section). In view of this, Theorem $1$ in this case
turns into Theorem $4.1$ in \cite{Gi0} and is a reformulation in terms of
the formalism explained in Section $3$ of the results of Mumford \cite{Mu} and 
Faber -- Pandharipande \cite{FP}. The proof of Theorem $1$ is based on a
similar application of Mumford's Grothendieck -- Riemann -- Roch argument to 
our somewhat more general situation. The argument, no doubt, was known 
to the authors of \cite{FP}. The main new observation here is that the
combinatorics of the resulting formula which appears rather complicated at
a first glance fits quite nicely with the formalism of quantized quadratic
hamiltonians. A verification of this --- somewhat tedious but 
straightforward --- is  presented in Appendix $1$. 

\bigskip

\nin {\bf 4. The Euler class.} 
The $S^1$-equivariant Euler class of $E$ is written in terms of the 
(non-equivariant) Chern roots $\rho_i$ as 
\[ \e (E) =\prod_i (\l +\rho_i). \]
Using the identity
$(\l+x)=\exp (  \ln \l - \sum_k (-x)^k / k \l ^k )$
we can express it via the components of the non-equivariant Chern character: 
\be \label{e} \e (E) = \exp \{ \ \ch_0(E) \ln \l + 
\sum_{k>0} \ch_k(E) \frac{(-1)^{k-1}(k-1)!}{\l^k}\ \} .\end{equation}   
Denote by $\D_{\e}$ the element $\D_{\s}$ of the Fock space corresponding
to $s_0=\ln \l$ and $s_k=(-1)^{k-1}(k-1)!/{\l^k}$ for $k>0$.
Substituting the values of $s_k$ into (\ref{T1}), 
replacing $\ch_l(E)$ by $\sum \rho_i^l/l!$ and using the binomial formula 
\[
(1+x)^{1-2m}=\sum_{l\geq 0} \frac{(-1)^l(2m-2+l)!}{(2m-2)!\ l!}\ x^l  
\]
we arrive at the following conclusion.

\medskip

{\bf Corollary 1.}
\[ \prod_i \exp \{ -\frac{1}{24} \int_X 
[(\l+\rho_i)\ln (\l+\rho_i)-(\l+\rho_i)]\ c_{D-1}(T_X) \} \ 
\prod_i (\operatorname{sdet} \sqrt{\l+\rho_i})^{-\frac{1}{24}}\ \D_{\e} =\]
\[ \prod_i \exp \{ 
\sum_{m>0} \frac{B_{2m}}{2m (2m-1)} \widehat{(\frac{z}{\l +\rho_i})^{2m-1}} \} 
\prod_i \exp \{ \frac{(\l+\rho_i)\ln (\l+\rho_i)-(\l+\rho_i)}{z} \}\hh
 \D_{\0} .\] 

\medskip

{\em Remark.} The $1/z$-term in this formula actually arises in the form 
\[ \rho \ln \l +\sum \frac{(-1)^{k-1}\rho^{k+1}}{k(k+1)\l^k} = 
\int_{0}^{\rho} \ln (\l+x) dx = [(\l+x) \ln (\l+x) - (\l+x)]\ | _{0}^{\rho} .\]
It has positive cohomological degree and is small in this sense.
The constant term $(\l \ln \l -\l)/z$ is thrown away on the following grounds. 
According to \cite{Gi0}, $(1/z)\hh $ is the {\em string operator}  
and annihilates the descendent potential $\D_{\0}$. 
Thus the operators $\exp ((\l \ln \l -\l)/z)\hh )$ do not change $\D_{\0}$.   
The rest of the series in the exponent converges in the $1/\l$-adic topology. 

\bigskip

\nin {\bf 5. Quantum Serre.} 
Introduce the multiplicative characteristic class
\[ \cc^*(\cdot ) = \exp\{ \sum (-1)^{k+1} s_k \ch_k (\cdot) \}.\]
Since $\ch_l(E^*)=(-1)^l\ch_l(E)$ we have
\[
\cc^*(E^*) = {1 \over \cc (E)} \ .
\]
There is no obvious relationship between $\cc^*((E^*)_{g,n,d})$ and
$\cc (E_{g,n,d})$, but nonetheless the twisted descendent potentials 
$\D_{\s}=\D_{\cc, E}$ and $\D^*_{\s}:=\D_{\cc^*, E^*}$ are closely related.

\medskip
\label{Cor2}
{\bf Corollary 2.} 
{\em We have $\D^*_{\s} = (\sdet{\cc (E)})^{-1/24} \D_{\s}$.
More explicitly, 
\[ \D_{\cc^*,E^*} (\t^* ) = \ (\sdet{\cc (E)})^{-\frac{1}{24}} 
\D_{\cc, E} (\t),
\ \text{where}\ \t^*(z) = \cc (E) \t (z) + (1-\cc (E)) z .\] }  

{\em Proof.}  Replacing $\ch_l(E)$ with $(-1)^l\ch_l(E)$, and $s_k$
with $(-1)^{k+1}s_k$ in (\ref{T1}) preserves all terms except the
super-determinant. 

\medskip

{\bf Corollary 3.} {\em Consider the dual bundle $E^*$ equipped
with the dual $S^1$-action, and the $S^1$-equivariant inverse Euler
class $\e^{-1}$.  Put
\[
\t^*(z) = z+(-1)^{\dim{E}/2} \e (E) (\t (z) -z)
\]
and introduce the change $\pm: Q^d \mapsto Q^d (-1)^{\int_d \ch_1(E)}$
in the Novikov ring. With this notation
\[ \D_{\e^{-1}, E^*} (\t^*, Q) = 
\sdet [(-1)^{\frac{\dim{E}}{2}}\e(E)]^{-\frac{1}{24}}
\D_{\e, E}(\t, \pm Q) .\]} 

\medskip

{\em Proof.} We have
$\e^{-1}(E^*)=\prod_i (-\l - \rho_i)^{-1}.$
Since
\[
(-\l+x)^{-1} =\exp \{ -\ln (-\l) +\sum_k {x^k \over k \l^k} \}
\]
we find that $\e^{-1}(\cdot) =\exp { \sum s_k^* \ch_k(\cdot)}$
where $s^*_k= (k-1)!/ \l^k$ for $k>0$ and $s_0=-\ln (-\l)$.
For $k>0$, $s_k^* = (-1)^{k+1} s_k$ as in the situation of Corollary
$2$.  However, $s^*_0= -s_0 -\pi \sqrt{-1}$.  We compensate for the
discrepancy $-\pi \sqrt{-1}$ using the divisor equation (\ref{de})
described in Remark $1$ following Theorem $1$.

\bigskip
 
\nin {\bf 6. Genus $0$.} The genus $0$ descendent potential $\F^0_X$ can be
recovered from the so called ``J-function'' of finitely many variables due
to a reconstruction theorem essentially due to Dubrovin \cite{Db} and
going back to Dijkraaf and Witten \cite{DW}. The {\em J-function}
is a formal function of $t\in H$ and $1/z$ with vector coefficients in $H$ 
defined by 
\be \label{J} \forall a\in H,\ \ 
(J_X (t, z), a):=(z+t,a)+\sum_{d,n} \frac{Q^d}{n!} \int_{[X_{0,n+1,d}]} 
\bigwedge_{i=1}^n \ev_i^*t \ \frac{\ev_{n+1}^*a}{z-\psi_{n+1} } . \end{equation}
We need the following reformulation of the reconstruction theorem
in terms of the geometry of the symplectic space $(\H, \O)$, where
we take $\H$ to be the completion $H((z^{-1}))$.

The genus-$0$ descendent potential $\F^0_X$ considered as
a formal function of $\q \in \H_{+}$ via the dilaton shift (\ref{ds})
generates (the germ of) a Lagrangian section $\L_X \subset \H = T^*\H_{+}$.
In Darboux coordinates 
\[ \L_X =\{ ( {\mathbf p}, \q ) : \ 
{\mathbf p }=d_\q \F^0_X \}\ .\] 

\medskip

{\bf Proposition.} {\em $\L_X$ is a homogeneous Lagrangian cone
swept by a moving semi-infinite isotropic subspace 
depending on $\dim H$ parameters. More precisely,  
 
(i) the tangent space $L_{\f} \subset \H$ to $\L_X$ at a point $\f$  
satisfies $\L_X \cap L_{\f} = z L_{\f}$;
 
(ii) $J_X(t,-z) \in \H$ is the intersection of $\L_X $ with 
$(t-z)+\H_{-}$.}

\medskip

{\em Remarks.} (1)  Part (i) implies that the tangent spaces $L_{\f}$ are 
Lagrangian subspaces invariant under multiplication by $z$. 
They consequently belong to the loop group Grassmannian 
(of the ``twisted'' series  $A^{(2)}$) or to its super-version.

 (2) Part (i) of the Proposition means that the spaces $L_{\f}$ actually
depend only on $\dim H$ parameters and form a {\em variation of semi-infinite
Hodge structures} in the sense of \cite{B}.    
It also shows that the cone $\L_X$ is determined by its generic 
$(\dim H)$-parametric slice $(J (t)\in \H)_{t\in H}$. Indeed, if the first
$t$-derivatives of $J$ span $L_{J(t)}/z L_{J(t)}$ over $\gL$, then they
span the tangent space $L_{J(t)}$ over $\gL [z]$ and the cone is the union 
of the isotropic spaces $z L_{J(t)}$. 
Part (ii) identifies one such slice with the J-function.
  
\medskip

Part (ii) of the Proposition follows immediately from the definitions 
of $J_X$ and $\L_X$. 

Part (i) follows easily from Dubrovin's reconstruction formula 
\cite{Db} in the axiomatic theory of Frobenius structures. 
Indeed, the main feature of the formula is that the
$2$nd differentials $d^2_{\q } \F_X^0$ of the genus $0$ descendent potential
depend on the application point $\q (z)=\t(z)-z$ only through some 
finite-dimensional function $\q \mapsto t(\q) \in H$, and the levels of the 
function are (germs at $\q=-z$ of) linear subspaces of codimension $\dim H$. 
In geometric terms this means that the tangent spaces to $\L_X$ (regarded as
affine spaces in $\H$) remain constant along these subspaces. The tangent spaces 
actually pass through the origin since $\L_X$ is a cone --- this
follows from the genus $0$  {\em dilaton equation} (see \cite{Db, Gi1}). 
The tangent spaces therefore form 
a family $\{ L_t \}$ of Lagrangian spaces depending only on 
$t\in H$ and intersecting the cone $\L_X$ along subspaces $I_t$ of 
codimension $\dim H$. Invariance of $\L_X$ with respect to the flow of the 
string vector field $\f \mapsto \f /z$ implies $z^{-1}I_t \subset L_t$. 
Considering the Fredholm index of projections to $\H_{+}$ shows that 
that $I_t=zL_t$ (since the codimension of $I_t$ in $L_t$ would exceed $\dim H$ otherwise). 

In the above argument, we assume that the ground ring $\gL$ is (or has been extended to) 
a field. In Appendix $2$ we give another, more direct proof  
applicable in Gromov -- Witten theory and free of this defect. 
It is based on Theorem $5.1$ stated in \cite{Gi0} which relates gravitational 
descendents with {\em ancestors}.

\medskip

In the quasi-classical limit $\h \to 0$,
quantized symplectic transformations $\exp \hat{A}$ of Theorem $1$ acting on 
the total potentials considered as elements $\D_{\s}$ in the Fock space turn 
into the ``unquantized'' symplectic transformations acting by 
$\L_{\s} \mapsto (\exp{A}) \L_{\s}$ on the Lagrangian 
cones $\L_{\s}$ generated by the genus $0$ potentials $\F^0_{\cc (\s), E}$. 

\medskip

{\bf Corollary 4.} 
\[ \L_{\s} = \exp \ \{ \ \sum_{m\geq 0}\ \sum_{0\leq l \leq D} s_{2m-1+l} 
\frac{B_{2m}}{(2m)!}  \ch_l(E) z^{2m-1}  \ \} \ \  \L_{\0} .\]

\bigskip

\nin {\bf 7. Quantum Lefschetz.} In the case of genus $0$ GW-theory
twisted by the Euler class $\e (E)$, the corresponding Lagrangian cone
$\L_{\s(\e)}$ is obtained from $\L_{X}$ by multiplication in $\H$ defined
by the product over the Chern roots $\rho$ of the series
\be \label{b} b_{\rho}(z)=
\exp \{ \frac{(\l+\rho )\ln (\l +\rho ) -(\l+\rho)}{z} + 
\sum_{m>0} \frac{B_{2m}}{2m(2m-1)} (\frac{z}{\l+\rho})^{2m-1} \} \end{equation} 
The series (\ref{b}) is well-known \cite{JEL} in connection with 
the asymptotic expansion of the gamma function $\Gamma ((\l+\rho)/z )$.
More precisely, (\ref{b}) coincides with the stationary phase 
asymptotics of the integral 
\[ 
\frac{1}{\sqrt{2\pi z (\l+\rho) }} \ \int_0^{\infty}
e^{\frac{-x+(\l+\rho)\ln x}{z}} dx \]
near the critical point $x=\l+\rho$ of the phase function.

Let us assume now that $E$ is the direct sum of $r$ line bundles with 
the $1$-st Chern classes $\rho_i$ --- in what follows we will need 
the Chern roots to be integer --- and consider the J-function 
$ J_X (t,z) = \sum_d J_d(t,z) Q^d $. Put $\rho_i(d)=\int_d\rho_i$ and
introduce the following {\em hypergeometric modification} of $J_X$:
\be \label{I} I_E(t,z) = \sum_d J_d (t,z) Q^d \prod_{i=1}^r
\frac{\prod_{k=-\infty}^{\rho_i(d)} (\l +\rho_i+k z)}{\prod_{k=-\infty}^{0}
(\l+\rho_i+kz)}. \end{equation} 
%Not accidentally, the infinite fractions here can be rewritten via the gamma function:
%$z^{\rho_i(d)} \Gamma (\rho_i(d)+1+(\l+\rho_i)/z) /\Gamma (1+(\l+\rho_i)/z)$. 
  
\medskip

{\bf Theorem 2.} {\em The hypergeometric modification $I_E$
considered as a family $t \mapsto I_E(t,-z)$ of vectors in the symplectic
space $(H, \O_{\e(E)})$ corresponding to the twisted inner product $(a,b)_{\e(E)}=
\int_X \e(E) ab$ on $H$, is situated on the Lagrangian section 
$\L_{\e,E} \subset \H$ defined by the differential of the twisted genus $0$ 
descendent potential $\F^0_{\e,E}$.}

\medskip

The following comment is in order. The series $I_E$ does not necessarily belong
to $H((z^{-1}))$ because of possible unbounded growth of the numbers
$\rho_i(d)$. However the coefficients at each particular monomial $Q^d$ do.
Similarly, multiplication by the series (\ref{b}) moves the cone $\L_X$ out
of the space $H((z^{-1}))$. However modulo each particular power of
$1/\l $ it does not (the invariance of the cone with respect to the string
flow $\exp (\l \ln \l - \l)/z$ is once again essential here). 
In fact all our formulas make sense as operations with generating 
functions, \this give
rise to legitimate operations with their coefficients, because of the 
presence of suitable auxiliary variables (such as $s_k$ in Corollary $3$,
$1/\l$ in (\ref{b}), $Q$ in (\ref{I}), etc.) More formally, this means that
(i) the ground ring $\gL$ in $H=H^*(X,\gL)$ should be completed in a suitable 
``adic'' topology, and 
(ii) in the role of the symplectic space $\H$,  we should 
take the space $H \{ z, z^{-1} \} $ of Laurent series $\sum_{k\in \ZZ} h_k z^k $
``convergent'' in the sense that $\lim_{k\to \infty} h_k \to 0$ in the topology of $\gL$.  
In the following proof we will have to similarly replace $\gL [z]$ by
$\gL \{ z \}$, and the ring $\gL $ should be also extended by  
$\sqrt{\l}$.

\bigskip

\nin {\bf 8. Proof of Theorem $2$.}
Due to the equivariance properties (see \cite{Gi1}, Section 6) 
of J-functions with respect to the string and divisor flows (\ref{de}) 
we have  
\[ J_X (t+\sum (\l+\rho_i) \ln x_i) = 
 e^{\frac{\sum (\l+\rho_i) \ln x_i}{z}} \sum_d J_d (t,z) Q^d 
\prod_i x_i^{ \rho_i (d)} .\]
Integrating by parts (as in the proof of the identity $\Gamma (x+1)=x \Gamma (x)$)
we find 
\be \label{int} (2\pi z)^{-\frac{r}{2}} \  
\int_0^{\infty}dx_1 ...\int_0^{\infty} dx_r\  
e^{-\sum x_i/z} J_X (t+\sum (\l+\rho_i) \ln x_i)  = \end{equation}
\[ I_E (t,z)\  \sqrt{\e (E)} \prod_i \frac{1}{\sqrt{2\pi z (\l+\rho_i)}}
\int_0^{\infty} 
e^\frac{-x_i+(\l+\rho_i)\ln x_i}{z} dx_i \ .\] 
We conclude that the asymptotic expansion of the integral (\ref{int})
coincides with $ I_E(t,z) \sqrt{\e (E)} \prod_i b_{\rho_i}(z)$.

The multiplication by $\sqrt{\e (E)}$ identifies the Lagrangian cone
$\L_{\e,E} \subset (\H, \O_{\e (E)})$ with its normalized incarnation 
$\L_{\s (\e)} \subset (\H,\O)$. Therefore Theorem $2$ is equivalent to the 
inclusion $I_E (t,-z) \sqrt{\e (E)} \in \L_{\s (\e)}$ and, due to Corollary 
$4$, --- to $I_E(t,-z)\sqrt{\e (E)} \prod b_{\rho_i}(-z) \in \L_{\0}=\L_X$. 
It remains to show therefore that the asymptotic expansion of the integral 
(\ref{int}) belongs to the cone determined by the J-function $J_X(t,z)$. 
In fact we will prove the following

\medskip
 
{\bf Lemma.} 
{ \em For each $t$, the asymptotic expansion of the integral
 (\ref{int}) differs from 
$\l^{{\dim E}/2} J_X(t^*,z)$ (at some other point
$t^*(t)$) by a linear combination of the first $t$-derivatives of $J_X$ at $t^*$ 
with coefficients in $z\gL \{ z \}$.}

\medskip 

For this, we are going to use another property of the J-function $J_X$
well-known in quantum cohomology theory and in the theory of Frobenius 
structures
(see for instance Section $6$ in \cite{Gi1} and \cite{Db}).
The first derivatives $\p_{\a} J_X$ satisfy 
the system of linear PDEs
\be \label{Dm} 
z \frac{\p}{\p t^{\a}}\frac{\p}{\p t^{\b}} J_X (t,z) = 
\sum_{\c} A_{\a\b}^{\c}(t) \frac{\p}{\p t^{\c}} J_X (t,z).\end{equation}
where we use a coordinate system $t=\sum t^{\a} \phi_{\a}$ on $H$. 
Indeed, following \cite{B} we can argue that the second derivatives are linear 
combination of the first derivatives over $z^{-1}\gL \{ z\}$ since  
infinitesimal variations of the
tangent spaces $L_{\f}$ spanned by $\p J_X/\p t^{\a}$ are to stay
inside $z^{-1}L_{\f}$, and that on the other hand the second derivatives
are contained in $\H_{-}$ since $J_X \in z+t + \H_{-}$. 
Further analysis reveals that
$A_{\a\b}^{\c}$ are structural constants of the quantum cohomology algebra
$\phi_{\a} \bullet \phi_{\b} = \sum A_{\a\b}^{\c} \phi_{\a}$. In particular,
$z \p_{1} J_X = J_X$ since $1\bullet = \operatorname{id}$
(we use here the notation $\p_{v}$ for the directional
derivative in the direction of $v\in H$ and take $v=1$).

We can interpret (\ref{Dm}) as the relations defining the D-module
generated by $J_X$, \this obtained from it by application of all differential
operators. Using
Taylor's formula $ J_X (t+  y \rho  ) = \exp ( y \p_{\rho} ) J_X (t)$  
we now view 
(\ref{int}) as the asymptotic expansion of the oscillating integral
taking values in this D-module:
\[  (2\pi z)^{-\frac{r}{2}} \  \int_0^{\infty}dx_1 ...\int_0^{\infty} dx_r\  
e^{\frac{-\sum x_i + \sum (\l+z\p_{\rho_i}) \ln x_i  }{z}} \ \ 
J_X (t , z) \ \sim \]
\be \label{exp}  
\prod_i  e^{ \frac{(\l+z\p_{\rho_i}) \ln (\l+z\p_{\rho_i}) -(\l +z\p_{\rho_i})}
{z} + \frac{1}{2}\ln (\l+z\p_{\rho_i}) + \sum_{m>0} \frac{B_{2m}}{2m(2m-1)} 
\frac{z^{2m-1}} {(\l+z\p_{\rho_i})^{2m-1}} } \  J_X (t, z) .\end{equation}
The exact form of the series in not relevant this time. What matters is that 
the relations (\ref{Dm}) in the D-module allow us to rewrite any high order 
derivation as a differential operator of first order and that composition of 
derivations coincides with the quantum cup-product $\bullet$ modulo higher order terms
in $z$:
\[ z\p_{v_1} ... z\p_{v_N} = z \p_{ v_1\bullet ... \bullet v_N} + o (z) ,\]
where $o(z)$ stands for a linear combination of $z\p_{\phi_{\a}}$ with
coefficients in $z\gL \{ z\}$.
Using this (and also the relation $\l J_X = z\p_{\l \cdot 1} J_X$ mentioned 
earlier) we see that \begin{eqnarray*}  
\prod_i e^{ \frac{(\l+z\p_{\rho_i}) \ln (\l+z\p_{\rho_i}) - 
(\l+z\p_{\rho_i}) }{z} } J_X(t,z) &=& 
\prod_i e^{ \partial_{[(\l+\rho_i\bullet)\ln (\l+ \rho_i\bullet )
-(\l+\rho_i\bullet)]1} + {o(z) \over z} } J_X(t,z)\\
&=& e^{\frac{o(z)}{z}} J_X(t^*,z) 
\end{eqnarray*}
where $t^*(t)=t + [\ \sum (\l+\rho_i\bullet)\ln (\l+ \rho_i\bullet ) 
-(\l+\rho_i\bullet)\ ] \ 1 .$ 

\nin Processing next the factor  
$e^{\frac{1}{2}\ln (\l+z\p_{\rho_i})}$, we take out $\sqrt{\l}$.
The remaining factor
$ e^{ \frac{1}{2} \ln (1+ z\p_{\rho_i}/\l)}$ together with the rest of the 
exponent in the asymptotic 
expansion (\ref{exp}) yields an expression of the type
$e^{o(z)/z} J_X(t^*,z) $ too. 
We conclude that the expansion (\ref{exp}) assumes the form
\[ \l^{\frac{\dim E}{2}}
J_X (t^*,z) + \sum_{\a} C_{\a} (t^*,z)\ z \p_{\phi_{\a}} J_X (t^*,z) ,\]
where the coefficients $C_{\a} (t^*, \cdot )$ are in $\gL \{ z\}$ as required.
                
\medskip

{\em Remark.} The proof of the Lemma actually shows that for any
phase function $\Phi (v)$ of $v\in H$ the asymptotics of the
oscillating integral $\int dv \  e^{\Phi (z\p_v ) /z} J_X(t,z)$ generates
the same cone as $J_X$.  
%It would be interesting to find out if this invariance property 
%of the genus $0$ potential $\F^0_X$ admits higher genus generalizations.

\medskip

{\bf Corollary 5.} {\em Let $\L_{\e,E} \subset (\H, \O_{\e(E)})$ be the Lagrangian 
cone determined by $I_E(t,-z)$ and let $L_t$ be the tangent space to $\L_{\e, E}$
at the point $I_E(t,-z)$. Then the intersection (unique due to some
transversality property) of $z L_t$ with the
affine subspace $-z+z \H_{-}$ coincides with the value $J_{\e,E}(\tau ,-z) \in
-z+\tau(t) + \H_{-}$ of the J-function corresponding to the $(\e,E)$-modified 
GW-theory. In other words, 
\be \label{j} J_{\e,E}(\tau, z) = 
I_E(t,z)+\sum_{\a} c_{\a}(t,z)\ z\ \p_{\phi_{\a}}I_E(t,z),
\ \text{where}\  c_{\a}(t, \cdot) \in \gL \{ z \}, \end{equation}
and $\tau (t) $ is determined by the asymptotics $z+\tau\ (\mod \H_{-})$  of 
the RHS.}

\medskip
 
{\em Remarks.} (1) The procedure of computing $J_{\e,E}$ in terms of $I_E$ is
reminiscent of the Birkhoff factorization $U(z,z^{-1})=V(z^{-1})W(z)$ 
in the theory of loop groups. Moreover, the procedure applied to 
the first derivatives of $I_{E}$ instead of $I_E$ {\em is} an example
of Birkhoff factorization. Indeed, the derivatives form a 
$\gL \{ z \}$-basis $U (z,z^{-1})$ in $L_t$, and $W(z)$ is the transition
matrix to another, canonical basis $V(z^{-1})\in 1 +\H_{-}$ formed by
the first derivatives of $J_{\e,E}$.  

(2) A by-product of Corollary $5$ is a geometrical description of
the ``mirror map'' $t\mapsto \tau $: the J-function obtained as the
intersection $L_t \cap (-z+z\H_{-})$ comes naturally parameterized by $t$
which may have little common with the projections $\tau-z$ of the intersection 
points along $\H_{-}$.   
 
\bigskip

\nin {\bf 9. Mirror formulas.} Let us assume now that the bundle $E$ (which
is still the sum of line bundles with first Chern classes $\rho_i$)
is {\em convex}, \this spanned fiberwise by global sections,
and apply the above results 
to the genus $0$ GW-theory for a complete intersection $j:Y\subset X$ 
defined by a global section.  
While the above proof of Theorem $2$ fails miserably in the  
limit $\l = 0$, the definition of the series $J_{\e,E}$ and $I_E$ and the 
relation between them described by Corollary $5$ survive the non-equivariant
specialization. Namely, at $\l=0$ the J-function $J_{\e,E}$ degenerates into
\[ J_{X,Y}(t,z) = z+t + \sum_{d,n}\frac{Q^d}{n!}\ (\ev_{n+1})_*
[ \frac{\e (E'_{0,n+1,d})}{z-\psi_{n+1}} \w_{i=1}^n\ev_i^*t ] ,\]
where $(\ev_{n+1})_*$ is the cohomological push-forward along the
evaluation map $\ev_{n+1}: X_{0,n+1,d}\to X$ and $\e$ is
the (non-equivariant!) Euler class.  Here $E'_{0,n+1,d}\subset E_{0,n+1,d}$ 
is the subbundle defined as the kernel of the 
evaluation map $E_{0,n+1,d} \to \ev_{n+1}^*E$ of sections (from 
$H^0(\S, f^*E)$) at the $n+1$-st marked point. \\
 
The function $J_{X,Y}$ is related to the GW-invariants of $Y$ by 
\be \e (E) J_{X,Y}(j^*u ,z) =_{H_2(Y)\to H_2(X)} 
j_*J_Y(u , z) \label{J_Y} \end{equation} 
since $[Y_{0,n+1,d}]=\e(E_{0,n+1,d})\cap [X_{0,n+1,d}]$ (see for instance
\cite{Kr}). The long subscript here is to remind us that the corresponding
homomorphism between Novikov rings should be applied to the RHS.

On the other hand, the series $I_E$ in the limit $\l=0$ 
specializes to
\[ I_{X,Y}(t,z)=\sum_d J_d(t,z) Q^d \prod_i \prod_{k=1}^{\rho_i(d)} 
(\rho_i+kz) \]
since $\rho_i(d)\geq 0$ for all degrees $d$ of holomorphic curves.
Passing to the limit $\l=0$ in Theorem $2$ and Corollary $5$ we obtain the 
following ``mirror theorem''.

\medskip

{\bf Corollary 6.} {\em The series $I_{X,Y}(t,-z)$ and 
$J_{X,Y}(\tau,-z)$ determine the same cone. In particular, 
the series $J_{X,Y}$ related to the J-function of $Y$ by (\ref{J_Y})
is recovered from $I_{X,Y}$ via the ``Birkhoff factorization procedure'' 
followed by the mirror map $t\mapsto \tau$ as described in Corollary $5$.}

\medskip

{\em Remark.} Corollary $6$ is more general than the
(otherwise similar) quantum Lefschetz hyperplane section theorems
by Bertram and Lee \cite{Ber, Lee} and Gathmann \cite{Gath}
for (i) it is applicable to arbitrary complete intersections $Y$ without
the restriction $c_1(Y)\geq 0$ and (ii) it describes the J-functions not only
over the {\em small} space of parameters $t\in H^{\leq 2}(X,\gL)$ but over
the entire Frobenius manifold $H^*(X,\gL)$. In fact the results of
\cite{Gath} allow one to deal with both generalizations and to compute
recursively the corresponding GW-invariants one at a time. What has been 
missing so far is the part that Birkhoff factorization plays in the
formulations.

\medskip    

Now restricting $J_{X,Y}$ and $I_{X,Y}$ to the small parameter space
$H^{\leq 2}(X,\gL)$ and assuming that $c_1(E)\leq c_1(X)$ we can derive the 
quantum Lefschetz theorems of \cite{BCKS, Kim, Ber, Lee, Gath}. A
dimensional argument shows that the series $I_{X,Y}$ on 
the small parameter space has the form 
\[ I_{X,Y}(t,z) = z F(t) +\sum G^i(t)\phi_i + O(z^{-1}) ,\]
where $\{ \phi_i \}$ is a basis in $H^{\leq 2}(X,\gL)$,  
$G^i$ and $F$ are scalar formal functions and $F$ is invertible 
(we have $F=1$ and $G^i=t^i$ when the Fano index is not too small).  
 
\medskip

{\bf Corollary 7.} {\em When $c_1(E)\leq c_1(X)$ the restriction of 
$J_{X,Y}$ to the small parameter space $\tau \in H^2(X,\gL)$      
is given by} 
\[ J_{X,Y}(\tau, z) = \frac{I_{X,Y}(t,z)}{F(t)},\ \text{where} \ 
\tau = \sum \frac{G^i(t)}{F(t)} \phi_i .\]

\medskip

The J-function of $X=\CC P^{n-1}$ restricted to the small parameter plane 
$t_0+t P$ (where $P$ is the hyperplane class generating the algebra 
$H^*(X,\gL)=\gL [P]/(P^n)$) takes on the form
\be \label{cp} J_{X} = z \ e^{(t_0+Pt)/z} \sum_{d\geq 0} \frac{Q^d e^{dt}}
{ \prod_{k=1}^{d} (P+kz)^n } .\end{equation}
For a hypersurface $Y$ of
degree $l$ in $\CC P^{n-1}$ we then have
\be \label{P} I_{X,Y} = z e^{(t_0+Pt)/z}\sum_{d\geq 0} 
 Q^d e^{dt} \frac{ \prod_{k=1}^{ld}(lP+kz)}
{ \prod_{k=1}^{d} (P+kz)^n } .\end{equation}

\medskip

{\bf Corollary 8.} {\em On the small parameter space 

(i) $J_{X,Y}(t_0,t,z)=I_{X,Y}(t_0,t,z)$ when $l<n-1$;

(ii) $J_{X,Y}(\tau_0,t,z)=I_{X,Y}(t_0,t,z), \ \tau_0 = t_0+l!Qe^t$, when 
$l=n-1$;

(iii) $J_{X,Y}(t_0,\tau ,z)=I_{X,Y}(t_0,t,z)/F(t), \ \tau = G(t)/F(t)$, 
when $l=n$, and the series $F$ and $G$ are found from the expansion
$I_{X,Y}= \exp (t_0/z) [ zF + G P  + O(z^{-1}) ]$.}

\medskip

Projecting $J_{X,Y}$ by $j^*$ onto the cohomology algebra 
$\gL [P]/(P^{n-1}) \subset H^*(Y,\gL)$ we recover the mirror theorem of 
\cite{Gi0}, and in the case $l=n=5$ --- the quintic mirror formula of 
Candelas {\em et al} \cite{COGP}.

\bigskip

\nin {\bf 10. Further comments.} 
 
{\em On quantum Riemann -- Roch.} The operators 
$\ch_l (E) z^{2m-1}$ commute. In the non-equivariant setting this property
is preserved under quantization for the operators with $m\geq 0$ which
occur in Theorem $1$. This is due to the nilpotency of $\ch_l(E)$ with $l>0$. 
Also, the summand with $l=1$ on the LHS of (\ref{T1}) is the only one left in 
this case. Thus formula (\ref{T1}) simplifies in the non-equivariant case:
\[   \D_{\s}  = 
 \( e^{s_0(c_1(E), c_{D-1}(T_X))} \sdet \sqrt{\cc (E)} \)^{\frac{1}{24}}\ 
%e^{s_0(\ch_1(E), c_{D-1}(T_X))} 
e^{ \sum_{m\geq 0}\sum_{l\geq 0} s_{2m-1+l} 
 \frac{B_{2m}}{(2m)!} \widehat{ \ch_l(E) z^{2m-1} } } \ \D_{\0} .\] 
The formula defines a formal group homomorphism from the group of 
invertible multiplicative characteristic classes to invertible operators
acting on elements of the Fock space. It would be interesting to find 
a quantum-mechanical interpretation of the normalizing factor in this formula. 
Since the Fock space should consist of top-degree forms on $\H_{+}$ 
rather than function, the super-determinant probably takes on the role 
of the Jacobian of our ``bare hands'' identification 
$\q \mapsto \sqrt{\cc (E)} \q$. We do not have however a plausible ``physical''
interpretation for the other factor.  

\medskip

{\em On the Lagrangian cones.}  
In the case of genus $0$ GW-theory of $X=pt$ the cone 
$\L_{pt}$ is generated by the family of functions in one variable $x$:
\[ F(x, \q):= \frac{1}{2}\int_0^x Q^2(u)\ du,\ \ \text{where}\ 
Q(x)=\sum q_k \frac{x^k}{k!} \ . \] 
In particular, under analytic continuation
the cone $\L_X$ acquires singularities studied in geometrical optics 
on manifolds with boundary (see for instance \cite{Ar, Gi4, Shc}) and
called {\em open swallowtails}. 
It would be interesting to study singularities of $L_X$ under analytic 
continuation and to understand significance of the relationship with 
geometrical optics.

According to some results and conjectures of \cite{Db} and \cite{Gi0}, 
the Lagrangian cones $\L $ corresponding to semisimple Frobenius structures 
are linearly isomorphic to a closure of the Cartesian products of $\dim H$ copies 
of $\L_{pt}$, and various models in genus $0$ 
GW-theory differ only by the position of the product with respect to the 
polarization. The same is true for the cones $\L_{\s}$ corresponding to the 
different twisted theories on the same $X$: according to Corollary $4$
they are obtained from each other by linear symplectic transformations. 

The transformations form the multiplicative group 
$\exp (\sum \tau_m z^{2m-1})$ 
where $\tau_m$ are even elements of the {\em algebra} $H$. 
The action of this group on the semi-infinite Grassmannian resembles the 
abstract grassmannian interpretation of the KdV hierarchy. 
It would be interesting to further this analogy.

\medskip

{\em On the mirror theory.} When $X=\operatorname{pt}$, the function
$J_{\operatorname{pt}} = \exp (\t_{0}/z)$. When $E=\CC^n$ is the trivial bundle
over the point, the integral (\ref{int}) turns into
\[ \int_0^{\infty} ... \int_0^{\infty} e^{-\frac{x_1+...+x_n}{z}} 
(x_1 ... x_n)^{\frac{\l}{z}} \ dx_1\w ... \w dx_n\ .\]
It would be interesting to find a ``quantum symplectic reduction theorem''
which would explain how this integral is related to the $J$-functions of toric 
manifolds $X$ (see \cite{Gi5}) obtained by symplectic reduction from $\CC^n$.
For example, when $X=\CC P^{n-1}=\CC^n // S^1$, components of the $J$-function 
(\ref{cp})   
coincide with the complex oscillating integral 
\be \label{cpi} J_{X} (t) = z \ \int_{\gamma \subset \{ u_1 ... u_n =e^t \} } 
e^{\frac{u_1+...+u_n}{z}} \frac{d\ln u_1\w ...\w d\ln u_n}{d t} \end{equation} 
over suitable cycles. For a degree $l\leq n$ hypersurface $Y\subset X$,
this yields integral representations for $I_{X,Y}$ and $I_Y$.
Indeed the $I$-function (\ref{P}) is proportional
to the convolution (\ref{int})
\[  \int_{0}^{\infty} dv\ e^{-\frac{v}{z}} J_X (t+l\ \ln v) = 
 \int_{ \{ u_1 ... u_n = v^l e^t \} }    
e^{\frac{u_1+...+u_n-v}{z}}  \frac{d v \w d\ln u_1 \w ...\w d\ln u_n}{d t} .\]
Using the change $u_i \mapsto u_i v$ for $i=1,...,l \leq n$, we transform it to
the ``mirror partner'' of $Y$:
\begin{multline*} \frac{1}{2\pi i} 
\int_{ \{ u_1 ... u_n =e^t \} } e^{(u_{l+1}+...+u_n)/z} 
\frac{d\ln u_1\w ... \w d\ln u_n}{(1-u_1-...-u_l)\ dt } = \\ 
\int_{ \{  u_1...u_n=e^t; \ u_1+...+u_l=1\} } 
e^{\frac{u_{l+1}+...+u_n}{z}}\frac{d\ln u_1\w ... \w d\ln u_n}
{d(1-u_1-...-u_l)\w dt}  .\end{multline*} 

Another question. 
According to physics literature \cite{W1}, the mirror maps $t\mapsto \tau$ 
arise from the mysterious {\em renormalization}. 
According to \cite{CK} the mathematical content of some important examples 
of renormalization in quantum field theory is Birkhoff factorization in
suitable infinite-dimensional groups. Are {\em renormalization}
and {\em Birkhoff factorization} synonymous?

\medskip

{\em On Serre duality.} 
%In a cohomological formulation, 
%the classical Serre duality between 
%$H^0(\Sigma, V\otimes K^{1/2})$ and
%$H^1(\Sigma , V^*\otimes K^{1/2})$ says that for a family $p: Y\to B$ 
%of curves $\ch_k (p_* (V\otimes K^{1/2})) = -(-1)^k \ch_k (p_* 
%(V \otimes K^{1/2}))$. 
In the genus $0$ theory, when $E$ is convex
and $E^*$ is concave, the sheaves $E_{0,n,d}$ and $-E^*_{g,n,d}$ are
vector bundles with the fibers $H^0(\S , f^*E)$ and $H^1(\S, f^*E^*)$ 
respectively. Using the Euler class of $E^*_{g,n,d}$ one obtains
Gromov --  Witten invariants of the non-compact total space $E^*X$ of the bundle
$E^*$. The invariants twisted with the Euler class of $E_{0,n,d}$ can be 
interpreted as genus $0$ Gromov -- Witten invariants of the ``super-manifold'' 
$(\Pi E)X$, \this the total space of the bundle $E$ with
the parity of the fibers reversed. 

The ``non-linear Serre duality'' phenomenon emerged 
in \cite{Gi1, Gi2} in the context of fixed point localization
for genus $0$ Gromov-Witten invariants of $(\Pi E)X$ and $E^*X$.
The duality was formulated as identification (modulo
minor adjustments such as $\l \mapsto -\l,\ Q\mapsto \pm Q$) of certain
genus $0$ potentials {\em written in Dubrovin's canonical coordinates} 
of the semi-simple Frobenius structures associated with the two 
modified theories.  According to \cite{Gi0,Gi3} the total 
descendent potential of a semi-simple Frobenius structure can be
described in terms of genus $0$ data presented in canonical coordinates. 
This implies a higher genus version of the quantum Serre duality principle 
whenever the fixed point
localization technique \cite{Gi3} applies. Corollary $2$ and
its particular case described by Corollary $3$ assert the 
principle in much greater generality and show that both
the localization technique and the reference to semi-simplicity and
canonical coordinates in this matter are redundant. 

%The additive correction $(1-\cc (E)) z$ 
%in the change $\t(z) \mapsto \t^*(z)$
%may have dramatic consequences. For example, in the proof of a toric mirror theorem 
%\cite{Gi2} based on the Serre duality principle,
%the correction claims the sole responsibility for the mirror maps. 

Theorem $2$, Corollary $5$ and the mirror formulas of Section $9$
have Serre-dual partners. 
Replacing $\e$ and $E$ with with $\e^{-1}$ and $E^*$ 
(equipped with the dual $S^1$-action as in
Section $5$) we should change the inner product to
$(a,\b)_{\e^{-1}(E^*)}=\int_X \prod (-\l -\rho_i)^{-1} ab =\int_X (-1)^r 
\e^{-1}(E) ab $ and $I_E$ --- to $I^*_{E^*}:=$
\[  \sum_d Q^dJ_d\prod_i\frac
{\prod_{k=-\infty}^0(-\l-\rho_i+kz)}
{\prod_{k=-\infty}^{-\rho_i(d)}(-\l-\rho_i+kz)} = \sum_d (\pm Q)^dJ_d \prod_i
\frac{ \prod_{k=-\infty}^{\rho_i(d)-1}(\l+\rho_i+kz)}
{\prod_{k=-\infty}^{-1}(\l+\rho_i+kz)}. \]  
When the classes $\rho_i$ are positive, the bundle $E^*$ is {\em concave}
in the sense that $H^0(\S, f^*E^*)=0$ for all compact curves $\S $ of any 
genus. The GW-invariants twisted by $(\e^{-1},E^*)$ admit 
the non-equivariant specialization $\l=0$ (in moduli
spaces $X_{g,m,d}$ of positive degrees $d\neq 0$). The reader can check
that the results of Section $9$, appropriately adjusted 
to the case of $I^*_{E^*}$, reproduce genus $0$ mirror results of \cite{Gi2}.

\section*{Appendix 1. The proof of Theorem 1}

An application of the 
Grothendieck -- Riemann -- Roch theorem to the bundle $\ev_{n+1}^*(E)$ 
over the universal
curve $\pi: X_{g,n+1,d} \to X_{g,n,d}$ yields the following equality
($r,l,a,b\geq 0$):
\begin{multline}
\ch_k(E_{g,n,d}) = \\ \pi_* \left[
\sum_{r+l=k+1} {B_r \over r!} \ch_l(\ev^*(E)) \cdot 
\( \psi^r -
\sum_{i=1}^n (\sigma_i)_* \psi_i^{r-1} + 
{1 \over 2} \iota_* \sum_{a+b=r-2} (-1)^a \psi_{+}^a \psi_{-}^{b} \) 
\right] \label{chk}
\end{multline}
Here $\psi=\psi_{n+1}$, $\sigma_i:X_{g,n,d} \to X_{g,n+1,d}$ is the
section of the universal family defined by the $i$-th marked point, 
$\iota $ is the embedding into $X_{g,n+1,d}$  of the stratum 
$X^{Sing}_{g,n,d}$ of virtual codimension $2$ formed by nodes
of the curves, and $\psi_{+},\psi_{-}$ denote the $1$-st Chern classes of the
line orbibundles over $X_{g,n,d}^{Sing}$ formed by the cotangent lines
to the two branches of the curves at the nodes.

In order to justify (\ref{chk}) let us first recall the 
Grothendieck -- Riemann -- Roch Theorem \cite{SGA, Fu}: 
\be \label{grr} \ch (p_* V) = p_* (\ch (V) \Td (\calt_{Y/B}) ) ,\end{equation}
where $V$ is a vector bundle on
$Y$ and  $p : Y \to B$ is a {\em local complete intersection morphism}. 
The latter hypothesis means (see \cite{Fu}) that for some 
(and hence for any) embedding $j: Y \subset M$ into a non-singular space, 
the embedding $j\times p: Y\subset M\times B$ has a normal bundle $N_Y$ 
(\this the normal sheaf is locally free).  The difference 
$j^*T_M\ominus N_Y$ then takes on the role of the virtual relative
tangent sheaf $\calt_{Y/B}$ which in fact does not depend on the choice 
of $j$. Under the hypothesis on $p$ there exists a complex $0\to A^0 \to ... \to A^N\to 0$
of locally free sheaves on $B$ with cohomological sheaves equal to $R^ip_*(V)$
(later we will explicitly describe such resolutions for $\pi_*E=E_{g,n,d}$).
The $K$-theoretic push forward $p_*V$ is defined as an element $\sum (-1)^i[A^i]$ in the 
Grothendieck group $K^0(B)$ of vector bundles on $B$, and $\ch (p_*V)$ denotes
the topological Chern character $\sum (-1)^i\ch (A^i) \in H^*(B;\QQ)$.

We will apply the theorem in the orbispace/orbibundle situation
which reduces to the following. The moduli orbispace $X_{g,n,d}$
(together with the universal stable map) 
can be described (see for instance \cite{FuP}) as the quotient $P/G$ of 
a space $P$ by a semisimple complex Lie group $G$ acting on $P$ algebraically with 
{\em at most finite stabilizers}. 
By definition, orbisheaves and orbibundles on $P/G$ are $G$-equivariant
sheaves and bundles on $P$. Their $G$-equivariant characteristic classes 
are elements of $H^*_G(P,\QQ )$ which coincides with $H^*(P/G, \QQ)$ 
since the action of $G$ on $P$ is almost free. Moreover, the $G$-space $P$
comes together with the $G$-equivariant universal family $C\to P$ of 
stable maps $\ev: C\to X$, and $\ev^*(E)$ is a $G$-equivariant bundle on $C$.
Equivariant sheaves and bundles induce ordinary
sheaves and bundles over finite-dimensional approximations to the 
homotopy quotients $C_G \to P_G$. The equivariant characteristic classes are
determined by the ordinary characteristic classes of such approximations.
Technically speaking, the formula (\ref{grr}) should be applied to 
these approximations. 
%(We refer to \cite{EG} for a systematic treatment
%of the equivariant Grothendieck -- Riemann -- Roch theorem.) 
The hypotheses needed in (\ref{grr}) are satisfied because (a) the projection 
$C\to P$ is a local complete intersection morphism and
(b) the space $C$ admits an equivariant
embedding into a non-singular space.

The statement (a) is a local property of the map and follows from
the fact that the universal family of curves $C\to P$ is flat and therefore
in a neighborhood of a nodal point of one of the fibers can be induced from 
the semi-universal unfolding  
\be \label{node} \CC^2\to \CC: (x,y)\mapsto xy \end{equation}
of the nodal singularity $xy=0$.

The statement (b) follows from the construction \cite{FP, FuP}
(in terms of Hilbert schemes) of an equivariant embedding of $C\to P$ 
into a larger flat family of curves $\tilde{C} \to \tilde{P}$ with
$\tilde{C}$ and $\tilde{P}$ non-singular. In particular, near the locus
$\tilde{C}^{Sing}$ of nodes the family of curves is transversally described
by the local normal form (\ref{node}). We will exploit this property
soon.

Further derivation of (\ref{chk}) does not differ much from Mumford's
argument \cite{Mu}. 
The normal form (\ref{node}) allows one 
to express  the relative cotangent orbisheaf 
$\calt^* = \calt^*_{X_{g,n+1,d}/X_{g,n,d}}$ 
via the universal cotangent line bundle $L$ at the last marked point. 
Put $\calo=\calo_{X_{g,n+1,d}}$. The sheaf $\calo (L)$ 
consists of meromorphic differentials on the curves 
allowed poles of order $\leq 1$ at the marked points and identified 
near the nodes with sections of the relative dualizing sheaf 
(which in the notations (\ref{node}) have the form 
$a(x,y)dx\w dy/d(xy)$). 
The sheaf $\calo (L)$ contains $\calt^*$ as a subsheaf of differentials 
holomorphic at the marked points and of the form 
$(b(x,y)x+c(x,y)y) dx\w dy/d(xy) \equiv bdy-cdx\ \mod \CC [[x,y]] d(xy)$ 
near the node (\ref{node}). 
The coordinate expression $dx\w dy/d(xy)$ represents a well-defined
locally constant section of the orbibundle 
$\CC_{\pm}:=
\Lambda^2(L_{+}\oplus L_{-})\otimes L_{+}^{-1}\otimes L_{-}^{-1}$ over the
singular locus $X_{g,n,d}^{Sing}$ where $L_{+}$ and $L_{-}$ are cotangent
lines at the nodes. Using this and the residue of meromorphic differentials 
at the marked points we find  
\[ \calo (L) / \calt^* = \iota_* \calo_{X_{g,n,d}^{Sing}}(\CC_{\pm}) \oplus 
(\sigma_1)_* \calo_{X_{g,n,d}} \oplus ... \oplus (\sigma_n)_* \calo_{X_{g,n,d}}
.\]
The class $c_1(L)=\psi $ vanishes when restricted to the pairwise disjoint 
strata $D_i=\sigma_i(X_{g,n,d})$ and $Z=\iota (X^{Sing}_{g,n,d})$. 
This translates the
multiplicative property of the dual Todd class $\Tdv (\cdot )$ to additivity
of $\Tdv (\cdot ) -1$:
\be \label{Tdv} \Td (\calt ) =\Tdv (\calt^*) = 1 + [\Tdv (\calo (L)) -1] 
+\sum_i [\frac{1}{\Tdv (\calo_{D_i})}-1] + [\frac{1}{\Tdv (\calo_Z)}-1 ] .\end{equation}
The first two terms yield 
\[  \Tdv (\calo (L) ) = \frac{\psi}{\exp \psi - 1} = \sum_{r\geq 0} \frac{B_r}
{r!}\psi^r .\]
Using $\sigma_i^*(-D_i)=\psi_i$ and the exact sequence 
$0\to \calo (-D_i) \to \calo \to \calo_{D_i} \to 0$ we find
\[ \frac{1}{\Tdv (\calo_{D_i})}-1 = \Td (\calo (-D_i)) -1 = \sum_{r\geq 1}
\frac{B_r}{r!} (-D_i)^r = - (\sigma_i)_* \sum_{r\geq 1} \frac{B_r}{r!}
\psi_i^{r-1}.\]
The codimension-$2$ summand in (\ref{Tdv}) is processed using the
inclusion-exclusion formula for the bi-graded Poincar\'e polynomial of 
$\CC [x,y]/(xy)$:
\[ \frac{1-uv}{(1-u)(1-v)}=\frac{1}{1-u}+\frac{1}{1-v}-1 .\]
Consider the enlarged nodal locus $\tilde{Z} \subset \tilde{C}/G$. On 
a double cover of its neighborhood $\tilde{Z}$ is the normal
crossing of the divisors $D_{\pm}$ with the conormal bundles $L_{\pm}$.
We see from the Koszul complex
\[ 0\to \calo (L_{+}\otimes L_{-}) \to \calo (L_{+})\oplus \calo (L_{-}) 
\to \calo \to \calo_Z \to 0 \]  
that in the neighborhood of $Z \subset X_{g,n+1,d}=C/G$
\[  {1\over \Tdv (\calo_Z)}-1 \ =\  \frac{1-e^{-D_{+}-D_{-}}}{D_{+}+D_{-}}
\frac{D_{+}}{1-e^{-D_{+}}}\frac{D_{-}}{1-e^{-D_{-}}} -1 \]
\begin{eqnarray*} 
% & = & \frac{1-e^{-D_{+}-D_{-}}}{D_{+}+D_{-}}
%\frac{D_{+}}{1-e^{-D_{+}}}\frac{D_{-}}{1-e^{-D_{-}}} -1 \\
 & = & \frac{D_{+}D_{-}}{D_{+}+D_{-}} ( {1\over 1-e^{-D_{+}}} + 
{1\over 1-e^{-D_{-}}} - 1 - {1\over D_{+}}-{1 \over D_{-}} ) \\
 & = & \frac{D_{+}D_{-}}{D_{+}+D_{-}} ( {1\over D_{+}}[\frac{D_{+}}{1-e^{-D_{+}}} -
\frac{D_{+}}{2}-1] +   {1\over D_{-}}[\frac{D_{-}}{1-e^{-D_{-}}} -
\frac{D_{-}}{2}-1] ) \\
 & = & \frac{1}{2}\iota_* [\sum_{r\geq 2} \frac{B_r}{r!} 
\frac{\psi_{+}^{r-1}+\psi_{-}^{r-1}}{\psi_{+}+\psi_{-}}] = \frac{1}{2}
 \iota_* [ \sum_{r\geq 2} \frac{B_r}{r!}
\sum_{a+b=r-2}(-1)^a\psi_{+}^a\psi_{-}^b ].\end{eqnarray*}
We use here  
$\psi_{\pm}=-\iota^*(D_{\pm})$, $B_0=1, B_1=-1/2,$, $B_r=0$ for odd 
$r>1$, and assume that the push-forward $\iota_*$ is taken with respect
to the virtual fundamental class $[Z]$ described in a neighborhood of $Z$ as
as the cap-product of $[X_{g,n+1,d}]$ with the Euler class $D_{+}D_{-}$ of 
the normal bundle of $\tilde{Z}$.
   
Combining the formulas for $\Td (\calt )$ with the Grothendieck -- 
Riemann -- Roch theorem we arrive at (\ref{chk}).\\

The formula (\ref{chk}) is the main geometric ingredient in out proof of 
Theorem $1$. We also need the following facts. 

(i) The comparison formula $\psi_i-\pi^*(\psi_i) = D_i$.

(ii) The naturality of the virtual fundamental cycles
$\pi^*[X_{g,n,d}]=[X_{g,n+1,d}]$ under the flat morphism $\pi $.

(iii) The composition rule for $X_{g,n,d}^{Sing}$. Namely, the singular locus 
coincides with the total range of the {\em gluing maps} 
\be \label{red}  
X_{g_{+},n_{+}+\bullet, d_{+}}\times_X X_{0,1+\bullet+\circ,0} \times_X
 X_{g_{-},n_{-}+\circ , d_{-}} \to 
X_{g,n,d}^{Sing} \subset X_{g,n+1,d} \end{equation}
(over all splittings $g=g_{+}+g_{-},\ n=n_{+}+n_{-},\ d=d_{+}+d_{-}$) and  
\be \label{irr} X_{g-1, n+\bullet+\circ} \times_{X\times X} 
X_{0,1+\bullet+\circ,0}  \to  X_{g,n,d}^{Sing} \subset X_{g,n+1,d} .\end{equation}
The composition rule says that images of the virtual fundamental classes
under the gluing maps add up to the virtual fundamental class $[Z]$ of
the singular locus.

These properties (ii) and (iii) are part of the axioms in \cite{KM}
proved in \cite{BF}, and (i) is well known too --- see for instance
\cite{W}.

Next, we need similar results about $E_{g,n,d}$ as elements
in the Grothendieck groups of coherent orbisheaves:
\begin{equation}
\tag{iv}	\pi^* E_{g,n,d} = E_{g,n+1,d}
\end{equation}
\begin{equation}
\tag{v}  \gamma^* \iota^* E_{g,n+1,d} = \operatorname{pr}_{+}^*
E_{g_{+},n_{+}+\bullet,d_{+}} +  
\operatorname{pr}_{-}^* E_{g_{-},n_{-}+\circ,d_{-}} - \ev_{\Delta}^* E 
\end{equation}
\begin{equation}
\tag{vi}	
\gamma^*\iota^* E_{g,n+1,d} =
E_{g-1,n+\bullet+\circ,d} - \ev_{\Delta}^*E 
\end{equation}
where $\gamma $ are the gluing maps (\ref{red}) and (\ref{irr}) respectively,
$\ev_{\Delta}=\ev_{\bullet}=\ev_{\circ}$ is the evaluation at the point
of gluing, and $\operatorname{pr}_{\pm}$ are projections to the factors.

The properties can be verified by representing the bundle $E$ on $X$ 
as the quotient $A/B$ of two concave bundles. For this, pick a 
positive line bundle $L$ and let the exact sequence 
$0\to Ker \to H^0(X; E\otimes L^N)\otimes L^{-N} \to E\to 0$ take on the role
of $0\to B\to A\to E\to 0$. Then $H^0(\gS ; f^*A)$ and $H^0(\gS ; f^*B)$ 
vanish for
sufficiently large $N$ and any non-constant $f: \gS \to X$ so that
$0\to H^0(\gS ; f^*E)\to H^1(\gS ; f^*B)\to H^1(\gS ; f^*A)\to H^1 (\gS ; f^*E)
 \to 0$ is exact. This construction applied to a universal stable map of degree $d\neq 0$
yields a locally free resolution $0\to R^1\pi_*(\ev^*B) \to R^1\pi_*(\ev^*A) \to 0$ for
$E_{g,n,d}=R^0\pi_*(\ev^* E) \ominus R^1\pi_*(\ev^*E)$ mentioned earlier   
and reduces the problem about the sheaves $E_{g,n,d}$ with $d\neq 0$ 
to the case of $H^1$-vector bundles $-A_{g,n,d}$ and $-B_{g,n,d}$.
When $d=0$, $H^0$ is non-zero but has constant rank too.    
In either case the 
formulas (iv), (v), (vi) are easy to check directly, for instance, using
Serre duality (it identifies elements of $H^1(\gS ; f^*A)^*$ with 
holomorphic differentials on $\gS-(\operatorname{nodes})$ with values in 
$f^*A$ and allowed 
poles of order $\leq 1$ at the nodes provided that the sum of the two
residues at each node equals zero).    

Finally, we will need three integrals over low-genus moduli spaces.
Let us introduce the following correlator notation:
suppose that
\[
\mathbf{a}^i(\psi) = a^{(i)}_0 + a^{(i)}_1 \psi + \ldots
\]
are polynomials in $\psi$ with coefficients in $H^*(X; \Lambda )$, and
$\beta \in H^*(X_{g,n,d}; \Lambda )$.  Write
\[
\< \beta ; \mathbf{a}^1, \ldots \mathbf{a}^n\>_{g,n,d} = 
\int_{[X_{g,n,d}]} \beta \w 
\(\sum_{j \geq 0} \ev_1^* (a_j^{(1)}) \ \psi_1^j \) 
\ \ldots \ 
\(\sum_{j \geq 0} \ev_n^* (a_j^{(n)}) \ \psi_n^j \)  .
\]
Also, set $ \cc_{g,n,d} = \exp \( \sum s_k \ch_k\(E_{g,n,d}\)\) \in
H^*(X_{g,n,d}; \Lambda )$.
In this notation
\be \tag{vii} \lan \cc_{0,3,0}; \t , \t, \ch_{k+1}(E) \ran_{0,3,0} = 
\int_X \cc (E)\ t_0^{\w 2} \ \ch_{k+1}(E)\ , \end{equation}
\be \tag{viii} \lan \cc_{1,1,0}; \ch_k(E) \psi \ran_{1,1,0} = 
 {1 \over 24} \int_X \e(X)\ \ch_k(E)\ , \end{equation}
\be \tag{ix} \lan \cc_{1,1,0}; \ch_{k+1}(E) \ran_{1,1,0} = \hspace{4cm} \end{equation} 
\[ \ \hspace{10mm} {1 \over 24} \int_X \e(X) \ch_{k+1}(E)
 \(\sum_j s_j \ch_{j-1}(E) \) - 
{1 \over 24} \int_X \cc (E) c_{D-1}(X) \ch_{k+1}(E) . \]
The equality (vii) is obvious since $[X_{0,3,0}]=[X]$, and (viii) and (ix) 
follow easily from the well-known facts:
$X_{1,1,0}=X\times \M_{1,1}$, $[X_{1,1,0}] =
\e (T_X \otimes \EE^{-1})\cap [X\times \M_{1,1}]$ (where $\EE $ is the Hodge
line bundle over the Deligne -- Mumford space $\M_{1,1}$),
$E_{1,1,0} = E \otimes (1 \ominus \EE^{-1})$, $c_1(\EE)=\psi$ and $\int_{[\M_{1,1}]} 
\psi = 1/24$.  \\
  
Using (\ref{chk}) and the properties (i -- ix) we now derive
Theorem 1. At $\s=0$, Theorem 1 holds trivially, so it suffices to prove the
infinitesimal version
\begin{multline}
{\partial \over \partial s_k} \D_{\s} = 
\(\sum_{\substack{2m+r=k+1 \\ m \geq 0}} {B_{2m} \over (2m)!} 
\widehat{(\ch_r(E) z^{2m-1})}\) \D_{\s} \\
+ \( \begin{array}{c}
{1 \over 24} \int_X c_{D-1}(X) \wedge \ch_{k+1}(E) 
+ {1 \over 48} \int_X \e(X) \wedge \ch_k(E) \\
- {1 \over 24} \int_X \e(X) \wedge \ch_{k+1}(E) \wedge 
\( \sum_l s_{l+1} \ch_l(E) \) 
\end{array} \) \D_{\s} \label{infl} \\
\end{multline}
Here the first two exceptional terms come from the factors on the LHS
of (\ref{T1}); in particular the second one is due to
\[
\(\sdet \sqrt{\cc(E)}\) = \exp\(\str \ln \sqrt{\cc(E)}\) = 
\exp\(\int_X \e(X) \wedge \({1 \over 2} \sum_j s_j \ch_j(E)\)\)
\]
The third exceptional term is the cocycle value
\[
\C ( \frac{B_2}{2}\sum s_{l+1} \widehat{\ch_l(E) z} , \widehat{\ch_{k+1}(E)/z} ) = 
-{1 \over 2} \str\(\ch_{k+1}(E) \cdot {1 \over 12} \sum_l s_{l+1} \ch_l(E) \)
\]
which arises from commuting the derivative of the $\frac{1}{z}$ terms
(on the RHS in (\ref{T1})) past the terms involving $z$. \\

In the above correlator notation, 
\[
\D_{\s} = \exp\(\sum_{g,n,d} \hbar^{g-1} { Q^d\over n!} 
\< \cc_{g,n,d} ; \t,\ldots,\t\>_{g,n,d}\)\ 
\]
and so
\begin{multline} \label{dD/ds}
\D_{\s}^{-1} {\partial \over \partial s_k} \D_{\s} = 
 \sum_{g,n,d} {Q^d \hbar^{g-1} \over n!} 
\< \ch_k(E_{g,n,d}) \wedge \cc_{g,n,d};\t,\ldots,\t\>_{g,n,d} \\
+  \sum_{g,n,d} {Q^d \hbar^{g-1} \over (n-1)!} 
\<\cc_{g,n,d};\t,\ldots,\t, {\partial \t \over \partial s_k}\>_{g,n,d}\ . 
\end{multline}
We apply our expression (\ref{chk}) for $\ch_k(E_{g,n,d})$ and
compare the result with (\ref{infl}) by extracting terms involving the same 
Bernoulli numbers. 

We begin with $B_0=1$. 
Due to the comparison  formula (i) we have 
\[ \pi^*(\t (\psi_i)) = \t (\psi_i) - (\sigma_i)_* [\t (\psi_i)/\psi_i]_{+},\]
where $[\cdot ]_{+}$ means power series truncation. Together with the 
naturality (iv) of the class $\cc_{g,n,d}$ under $\pi^*$, this implies that 
\begin{eqnarray*} 
\lan \pi_*\left[ \ev^*\ch_{k+1}(E)\right] \cc_{g,n,d}; \t,...,\t \ran_{g,n,d} 
& = & \\
\lan \cc_{g,n+1,d}; \t, ..., \t, \ch_{k+1}(E) \ran_{g,n+1,d} 
& - & n \lan \cc_{g,n,d}; \t, ...,\t, \ch_{k+1}(E) 
\left[ \frac{\t(\psi)}{\psi}\right]_{+} \ran_{g,n,d} .\end{eqnarray*}
Summing over $g,n,d$ we find
 \begin{multline}
 \sum_{g,n,d} {\hbar^{g-1} Q^d \over n!} 
\< \pi_*\left[ B_0\ev^*\ch_{k+1}(E) \right] \cc_{g,n,d};
\t,\ldots,\t \>_{g,n,d} 
 = \\  - \sum_{g,n,d} {\hbar^{g-1} Q^d\over (n-1)!} 
 \< \cc_{g,n,d};\t,\ldots,\t, \ch_{k+1}(E) 
\left[ {\t(\psi)-\psi \over \psi} \right]_{+} \>_{g,n,d} \\ 
 - \frac{1}{2\h } \< \cc_{0,3,0};\t,\t,\ch_{k+1}(E)\>_{0,3,0} 
 - \<\cc_{1,1,0};\ch_{k+1}(E)\>_{1,1,0} \ .
 	  \label{B0detail}
\end{multline}
Here the exceptional terms arise from the fact that the moduli spaces
$X_{0,2,0}$ and $X_{1,0,0}$ are empty and therefore
$X_{0,3,0}$ and $X_{1,1,0}$ cannot be interpreted as universal curves.

The first two summands on the right actually add up to  
$\D_{\s}^{-1} (\ch_{k+1}/z)\hh \D_{\s}$. Indeed, the corresponding quadratic 
hamiltonian has $pq$-, $q^2$-, but no $p^2$-terms. Quantization of the 
$pq$-terms yields a linear vector field defined by the operator
$\q (z) \mapsto -[ \ch_{k+1}(E) \q (z)/z ]_{+}$, while the
$q^2$-term is $-(q_0,q_0)/2$ and matches the $2$-nd summand in (\ref{B0detail})
due to (vii). Evaluating the third summand
via (ix) we conclude that (\ref{B0detail}) coincides with
\be
\D_{\s}^{-1} \left( \begin{array}{l}
\widehat{\({\ch_{k+1}(E) \over z}\)} 
+ {1 \over 24} \int_X c_{D-1}(X)\ \ch_{k+1}(E) \\
- {1 \over 24} \int_X \e(X)\  \ch_{k+1}(E) \   
\( \sum_j s_j \ch_{j-1}(E) \) 
	\end{array} 
\right) \D_{\s} \ .
\label{B0}
\end{equation} 

Next, proceeding as above with $B_1=-1/2$ 
and using $\sigma_i^*\psi =0$ we find 
\begin{multline*}
 \sum_{g,n,d} {\hbar^{g-1}Q^d \over n!} 
\lan \pi_* \left[ B_1 \ev^* \ch_k(E)  \( \psi - D_1 - \ldots - D_n \)
\right] \cc_{g,n,d} ; \t, \ldots, \t \ran_{g,n,d} = \\
 \frac{1}{2} \sum_{g,n,d} {\hbar^{g-1}Q^d \over (n-1)!} \lan \cc_{g,n,d};
\t, \ldots, \t, \ch_k(E) ( \t (\psi ) - \psi ) \ran_{g,n,d}  
 +\frac{1}{2} \lan \cc_{1,1,0} ; \ch_k(E) \psi \ran_{1,1,0} .
\label{B1detail} 
\end{multline*}
In view of (viii) this coincides with
\be   - \sum_{g,n,d}
{ \hbar^{g-1} Q^d \over (n-1)!}  
\<\cc_{g,n,d};\t,\ldots,\t, {\partial \t \over \partial s_k}\>_{g,n,d} 
+ {1 \over 48} \int_X \e(X) \ \ch_k(E)  
\label{B1} \end{equation}
since $\t(z) = \cc (E)^{-1/2} \q(z) + z$ and hence $\p \t (z) /\p s_k = 
-\ch_k (E) (\t (z)-z)/2$.

Finally, it remains to check the equality of the $B_{2m}$-terms with $m>0$: 
\begin{multline} \label{B2m}
 \sum_{g,n,d} {\hbar^{g-1} Q^d \over n!} 
\< \pi_*\left[ \ev^*\ch_{k+1-2m}(E) \ 
\Psi_m \right] \ \cc_{g,n,d}; \t,\ldots,\t \>_{g,n,d} \\ 
=  \D_{\s}^{-1} \widehat{\(\ch_{k+1-2m}(E) z^{2m-1}\)} \D_{\s} 
\end{multline}
where 
\[
\Psi_m =  \psi^{2m} - \sum_{i=1}^n(\sigma_i)_* \psi_i^{2m-1} + {1
\over 2} \iota_* \( \frac{\psi_{+}^{2m-1}+\psi_{-}^{2m-1}}{\psi_{+}+\psi_{-}}
\) . \]
Processing the first two summands in $\Psi_m$ as before yields
\[  -\sum_{g,n,d} {\hbar^{g-1} Q^d \over (n-1)!}
\<\cc_{g,n,d}; \t,\ldots,\t, \ev^*\ch_{k+1-2m}(E) \ \psi^{2m-1} \ 
(\t(\psi) - \psi) \>_{g,n,d} ,\]
which coincides with the derivative of $\ln \D_s$ along the linear vector
field defined by the multiplication operator 
$\q (z) \mapsto -\ch_{k+1-2m}(E) z^{2m-1} \q(z) $. This vector field is
the quantization of the $pq$-terms in the quadratic hamiltonian 
corresponding to $\ch_{k+1-2m}(E) z^{2m-1}$. For $m>0$ 
the hamiltonian contains no $q^2$-terms. 

Let us identify the cohomology space $H$ with its dual by means of the
intersection pairing $(\cdot ,\cdot )$. We 
use the coordinate notation $\sum_{\a\b} \p_{\a} \ch^{\a\b} \p_{\b}$ for
the bi-derivation on $H$ corresponding to the self-adjoint operator 
of multiplication by $\ch_{k+1-2m}(E)$ on $H^*=H$. 
Applying the composition rules (iii),(v),(vi) we can express
contributions of the last summand in $\Psi_m$ as
\[ \D_{\s}^{-1} \left[ \frac{\h }{2} \sum_{a+b=2m-2}(-1)^a \sum_{\a\b} 
\p_{q_a^{\a}} \ch^{\a\b} \p_{q_b^{\b}} \right] \D_{\s} .\]
(In particular, the factor $\cc(E)$ due to $\p_t = \sqrt{\cc(E)} \p_q$
cancels with $1/\cc (\ev_{\Delta}^*(E))$ from (v) and (vi).)
This matches up with the quantization of $p^2$-terms in the quadratic 
hamiltonian of $\ch_{k+1-2m}(E)z^{2m-1}$. \\

Combining (\ref{B0}), (\ref{B1}), (\ref{B2m}) with (\ref{dD/ds}) and 
(\ref{infl}) we recover Theorem $1$.
   
\section*{Appendix 2. Descendents and ancestors}

The aim of the appendix is to justify part (i) of the Proposition in Section $6$ 
describing properties of the genus $0$ descendent potential $\F_X^0$ in 
terms of the geometry of the symplectic space $(\H,\O)$. In fact we intend to do 
more, namely --- to derive the Proposition from a relationship between 
gravitational descendents of any genus and their counterparts from  
Deligne -- Mumford spaces --- {\em ancestors} ---  expressed
in terms of the quantization formalism of Section $2$. 
The theorem in question, which is a reformulation of a result by 
Kontsevich -- Manin \cite{KM}, has been announced \cite{Gi0}. 
We recall the formulation and furnish a proof below.   

\medskip

Consider the composition $X_{g,m+l,d} \to \M_{g,m}$ of the 
operations of forgetting
the last $l$ marked points and contraction. Denote by $\bar{\psi}_i$ the pull-backs from 
Deligne -- Mumford space $\M_{g,m}$ of
the $1$-st Chern classes of universal cotangent lines. They differ from
the descendent classes $\psi_i$ on $X_{g,m+l,d}$. Following \cite{Gi0},
introduce the genus $g$ {\em ancestor} potentials
\be \label{a^g} \bar{\F}^g_X:= \sum_{d,m, l}^{\ }\frac{Q^d}{m! l!}
\int_{[X_{g,m+l,d}]} \bigwedge_{i=1}^ m 
[\ \sum_{k\geq 0}(\ev_i^*\bar{t}_k)\ \bar{\psi}_i^k\ ]\ 
 \bigwedge_{i=m+1}^{m+l}\ev_i^*\tau \ , \end{equation}
which are formal functions of the ancestor variables $\bar{\t}=
\sum \bar{t}_k\bar{\psi}^k$, $\bar{t}_k \in H$, and of the parameters 
$\tau\in H$.
The {\em total ancestor potential} is defined as
\[ \A_{\tau} = \exp \{ \sum \h^{g-1} \bar{\F}^g_X \} \] 
and is identified via the dilaton shift $\q (z) = \bar{\t} (z) - z$ 
with an element in the Fock space depending on the parameter $\tau\in H$.

We will use the abbreviated correlator notation
\[ \lan {\mathbf a}_1(\psi,\bar{\psi}) , ..., {\mathbf a}_m(\psi,\bar{\psi}) 
\ran_{g,m} (\tau ) := \sum_{l,d} \frac{Q^d}{l!}
\lan 1; {\mathbf a}_1(\psi,\bar{\psi}) , ..., {\mathbf a}_m(\psi,\bar{\psi}), 
 \tau, ... , \tau \ran_{g,m+l,d} \]
for Taylor series in $\tau $ with coefficients possibly mixing 
descendent and ancestor classes.
 
Introduce the operator series $S_{\tau}(z^{-1})=1+S_1z^{-1}+S_2z^{-2}+...$
acting on the space $\H = H ((z^{-1}))$ 
and defined in terms of genus $0$ descendents: 
\be \label{S} (S_{\tau} (z^{-1}) u, v):=
(u,v)+\lan \frac{u}{z-\psi}, v\ran_{0,2}(\tau ) .\end{equation} 
The series $S_{\tau}$ depends on the parameter $\tau\in H$. According to \cite{Gi1, Gi2}
it satisfies the identity $S^*_{\tau}(-z^{-1})S_{\tau}(z^{-1})=1$
 and consequently defines a symplectic transformation on 
$(\H, \O)$.
By quantization $\hat{S}$ of symplectic transformations we mean 
$\exp \widehat{\ln S} $. 

The action of the operator 
$\hat{S}_t^{-1}$ on an element $\G$ of the Fock space is explicitly described
by the formula:
\[ (\hat{S}_{\tau}^{-1} \G ) (\q) = e^{\lan \q,\q \ran_{0,2}(\tau)/2\h} 
\G ([S_{\tau}\q ]_{+}) ,\]
where $[S_{\tau}\q]_{+}$ is the power series truncation of $S_{\tau}
(z^{-1})\q(z)$. The formula is easy to check by generalizing it to
$\exp (-\epsilon \hat{A}) $ as in Proposition $5.3$ in \cite{Gi0} 
and taking $A=\ln S$, $\epsilon =1$. 
The quadratic hamiltonian of $A$ contains no $p^2$-terms 
(since $S$ is a power series in $1/z$), and quantization of $\exp (- \epsilon A)$ 
amounts to solving a $1$-st order linear PDE by the method of characteristics. 
The $pq$-terms give rise to the linear change of variables $\q \mapsto 
 [\exp (\epsilon A) \q]_{+}$. The exponential  factor can be verified ---
by differentiation in $\epsilon$  --- using the WDVV-like identity
\[ \lan \q (\psi) , 1, \q (\psi) \ran_{0,3}(\tau) = 
\sum_{\a\b} \lan \q (\psi), 1, \phi_{\a}\ran_{0,3}(\tau) g^{\a\b} 
\lan \phi_{\b}, 1, \q (\psi ) \ran_{0,3}(\tau) \]
(where $(g^{\a\b})$ is the inverse to the intersection matrix $g_{\a\b}=
(\phi_{\a},\phi_{\b})$) together with the string equation. We leave
some details here to the reader.

Let $F^1(\tau) :=\lan\ \ran_{1,0}(\tau)=
\F^1_X(\t )|\ _{ t_0=\tau , t_1=t_2=...=0}$ denote
the genus $1$ (non - descendent) GW-potential of $X$. Recall that the descendent 
potential $\D=\D_{X}$ is identified with an element of the Fock space
via the dilaton shift $\q(z) =\t (z) -z$.

\medskip

{\bf Theorem.} $\D = e^{F^1(t)} \hat{S}_{\tau}^{-1} \A_{\tau}$.

\medskip

{\em Proof.} Let $L$ be a universal cotangent line bundle over 
$X_{g,m+l,d}$, and $\bar{L}$ be its counterpart pulled back from $\M_{g,m}$ 
and corresponding to the same index (let it be $1$) of the marked point.
Let $\psi=c_1(L)$ and $\bar{\psi}=c_1(\bar{L})$.
There exists a section of $\mathop{Hom} (\bar{L},L)$ regular outside some virtual
divisor $D$ consisting of stable maps with the following property: 
the $1$-st marked point $\1$
is situated on a component of the curve which is subject to contraction 
under the map $X_{g,m+l,d}\to\M_{g,m}$.
It is easy to see that $D$ is the total range of the gluing maps
\[ X_{0,\1+\bullet+l',d'} \times_{X} X_{g,m-\1+\circ+l'', d''} \to X_{g,m+l,d} 
\] 
over all splittings $l'+l''=l, d'+d''=d$.
The virtual normal bundle to $D$ (outside self-intersections of $D$)
is canonically identified with $\mathop{Hom}(\bar{L},L)$. This implies 
that the section vanishes on $D$ with $1$-st order, and hence that 
$\psi-\bar{\psi}$
is Poincar\'e-dual to the virtual divisor: $[X_{g,m+l,d}]\cap (\psi-\bar{\psi})
=[D]$. Thus we have 
\begin{multline*} \lan u \psi^{a+1}\bar{\psi}^b, ... \ran_{g,m}(\tau )\ = \\ 
\lan u \psi^{a}\bar{\psi}^{b+1}, ...\ran_{g,m}(\tau)\ +\  \sum_{\a\b}
\lan u \psi^a,\phi_{\a}\ran_{0,2}(\tau)\ g^{\a\b}\ \lan \phi_{\b} 
\bar{\psi}^b, ... \ran_{g,m}(\tau) , \end{multline*}
where dots mean the descendent/ancestor content of other marked points 
(to be the same in all three places).
Applying this identity inductively to reduce 
descendents to ancestors we conclude that the descendent potentials
$\lan \t(\psi),...,\t(\psi)\ran_{g,m}(\tau)$  are 
obtained from the corresponding ancestor potentials 
$\lan \bar{\t}(\bar{\psi}),...,\bar{\t}(\bar{\psi})\ran_{g,m}(\tau)$ 
by the substitution $\bar{\t}(z) = [S_{\tau}(z^{-1}) \t(z)]_{+} $. 
This is essentially the result from \cite{KM}.

Let us compare this conclusion with the statement of the theorem.
Noting the presence of the similar change $\q \mapsto [S_{\tau}\q ]_{+}$ in the
explicit description of the operator $\hat{S}_{\tau}^{-1}$ we should also 
notice that $\q$ and $\t $ are not the same: $\q(z)=\t(z)-z$. 
This gives rise to the discrepancy 
$[z-S_{\tau} z]_{+}$. Expanding 
\[ [S_{\tau}]^{\a}_{\b}=\d^{\a}_{\b}+z^{-1} \sum_{\mu}g^{\a\mu}
\lan \phi_{\mu},\phi_{\b}\ran_{0,2}(\tau)+ o (z^{-1}) ,\]
we find the discrepancy equal to $-\tau $ because in components  
\[  \sum_{\mu} g^{\a\mu} \lan \phi_{\mu}, 1\ran_{0,2} (\tau ) = 
\sum_{\mu}g^{\a\mu}\lan 1; \phi_{\mu},1,\tau \ran_{0,3,0}=
\sum_{\mu}g^{\a\mu}g_{\mu\b}\tau ^{\b}=\tau ^{\a}. \]
Thus $\bar{\q}=S_{\tau}\q$ is equivalent to 
\[ \bar{\t} = [S_{\tau}\t]_{+} -\tau = [S_{\tau}(\t - \tau)]_{+} .\]
By Taylor's formula, we have 
\[ \F^g (\t )= \sum_{m=0}^{\infty} \frac{1}{m!}\lan \t(\psi ),...,\t(\psi )\ran_{g,m} (0)
= \sum_{m=0}^{\infty} \frac{1}{m!}\lan \t(\psi )-\tau,...,\t(\psi)-\tau \ran_{g,m}(\tau) .\] 
We conclude that for $g>1$ the descendent potentials $\F^g$ 
(which do not depend on $\tau$) are obtained from the ancestor potentials 
$\bar{\F}^g$ (which do depend on $\tau$) by the substitution 
$ \bar{\q}(z)=[S_{\tau }(z^{-1})\q(z)]_{+}$. In order to make the same true for $g=0,1$ we
have to include the terms corresponding to the unstable indices 
$(g,m)=(0,0),(0,1),(0,2)$ and $(1,0)$ and hence missing from the ancestor potentials. 
The first three of them give rise to the factor 
$\exp \lan \q , \q \ran_{0,2}(\tau) /2\h$. Indeed, 
\[ \frac{1}{2}\lan \t(\psi)-\psi,\t(\psi)-\psi\ran_{0,2} (\tau) = 
\lan\ \ran_{0,0}(\tau )+\lan \t(\psi )-\tau \ran_{0,1}(\tau )+
\frac{1}{2}\lan\t(\psi )-\tau,\t(\psi )-\tau\ran_{0,2}(\tau) .\]
This can be easily derived from the {\em dilaton equation} 
$\lan \psi , ... \ran_{g,n+1,d}=(2g-2+n)\lan ... \ran_{g,n,d}$    
applied with $g=0$.  Finally, the missing summand $\lan \ \ran_{1,0}(\tau)$ 
coincides with $F^1(\tau)$, and the Theorem follows.

\medskip

Passing to the quasi-classical limit $\h \to 0$ we obtain the following
result.

\medskip

{\bf Corollary.} {\em The Lagrangian sections $\L$ and $\bar{\L}_{\tau}$ 
which represent respectively the differentials of the genus $0$ 
descendent potential
$\F^0$ and ancestor potentials $\bar{\F}^0_{\tau}$ are related by
the symplectic transformations: $ \bar{\L}_{\tau} = S_{\tau} \L$.}

\medskip

Finally, we derive part (i) of the Proposition. 

When the ancestor variable $\bar{\q}$ belongs to $z\H_{+}$ (\this 
$\bar{t}_0=0$), the genus $0$ ancestor potential $\bar{\F}^0$ has 
identically zero $2$-jet at $\bar{\q}$. This follows from 
$\dim \M_{0,m+2} < m$. Thus the cone $\bar{\L}_{\tau}$: (a) contains the 
isotropic space $z \H_{+}$ and (b) at any point $\bar{\q}\in z\H_{+}$ has
the tangent space $\bar{L}_{\bar{\q}} = \H_{+}$. Applying the 
symplectic transformation $S^{-1}_{\tau}$ we see that the tangent spaces 
$L_{\f}$ to $\L$ at $\f= S^{-1}_{\tau} \bar{\q}$ intersect $\L$ along 
$zL_{\f}$ provided that $\bar{\q} \in z \H_{+}$. The condition 
$S_{\tau}\f \in z \H_{+}$ on $\f=({\mathbf p},\q) \in \L$
is equivalent to the system of equations 
\[ \lan 1, \q (\psi), v \ran_{0,3} (\tau) =0 \ 
\text{for all} \ v \in H . \]
In other words, $\tau$ must be a critical point of 
$ \lan 1, \q (\psi) \ran_{0,2} (\tau)$ considered as a function of $\tau \in H$
(depending on the parameter $\q \in \H_{+}$).
When $\q(z)=q_0-z$, the function 
turns into $(q_0,\tau)-(\tau,\tau)/2$ and has the nondegenerate critical point
$\tau = q_0$. This guarantees existence of a unique critical point 
$\tau (\q)$ in a formal neighborhood of $\q =-z$. The result follows.

\enddocument